\documentclass[11pt]{article}

\usepackage{amssymb}
\usepackage[english]{babel}
\usepackage{mathrsfs}
\usepackage{extarrows}
\usepackage{amsthm}
\usepackage{setspace}
\usepackage{calligra}

\usepackage{amssymb, amsmath, amsthm}
\usepackage{bbm}
\usepackage{graphicx,verbatim}
\usepackage{rotating}
\usepackage[all,cmtip]{xy}
\usepackage[a4paper,left=1in,width=458pt,textheight=650pt,top=1.5in]{geometry}
\usepackage{tikz}
\usetikzlibrary{calc}
\usepgflibrary{shapes.geometric}
\usepgflibrary{shapes.misc}
\usetikzlibrary{positioning}
\usetikzlibrary{decorations}
\usetikzlibrary{arrows}

    \theoremstyle{plain}
    
    \newtheorem{theorem}{Theorem}[section]
    \newtheorem{proposition}[theorem]{Proposition}
    \newtheorem{corollary}[theorem]{Corollary}
    \newtheorem{lemma}[theorem]{Lemma}
    \newtheorem{conjecture}[theorem]{Conjecture}
    \newtheorem*{theorem*}{Theorem}
    \newtheorem*{proposition*}{Proposition}
    \newtheorem*{corollary*}{Corollary}
    \newtheorem*{lemma*}{Lemma}
    \newtheorem*{conjecture*}{Conjecture}

    \theoremstyle{definition}
    \newtheorem{definition}[theorem]{Definition}
    \newtheorem{example}[theorem]{Example}
    \newtheorem*{definition*}{Definition}
    \newtheorem*{example*}{Example}

    \theoremstyle{remark}
    \newtheorem{remark}[theorem]{Remark}
    \newtheorem*{remark*}{Remark}

    \newenvironment{myproof}[1][\noindent\textit{\sc \textbf{Proof.} }]
    {#1}{$\hfill\qedsymbol\newline$}

    \newcommand{\bD}{\mathbb{D}}
     
    \newcommand{\bZ}{\mathbb{Z}}

    \newcommand{\bA}{\mathbb{A}}

    \newcommand{\bK}{\mathbb{K}}
    \newcommand{\cC}{\mathcal{C}}
    \newcommand{\cR}{\mathcal{R}}
    \newcommand{\cO}{\mathcal{O}}
    
    \newcommand{\cD}{\mathcal{D}}

    \newcommand{\cU}{\mathcal{U}}

    \newcommand{\fM}{\mathfrak{M}}

    \newcommand{\fg}{\mathfrak{g}}
    \newcommand{\fh}{\mathfrak{h}}
    
    \newcommand{\fk}{\mathfrak{k}}
    \newcommand{\ft}{\mathfrak{t}}
    
    \newcommand{\sO}{\mathscr{O}}	
    \newcommand{\sD}{\mathscr{D}}
    \newcommand{\sA}{\mathscr{A}}
    \newcommand{\sV}{\mathscr{V}}
    \newcommand{\sY}{\mathscr{Y}}

    \renewcommand{\ker}{\operatorname{ker }}

     \newcommand{\qdiff}{\mathscr{D}_{q}}
     \newcommand{\qdiffc}{\mathscr{D}^{\circ}_{q}}
     \newcommand{\adiff}{\mathscr{A}^{\chi}_{q,\eta}}
     \newcommand{\ml}{\mathfrak{M}^{\paren{l}}_{\chi,\eta}}
     \newcommand{\mlaff}{\Phi_{K^{(l)}}^{-1}\paren{\eta^{l}}\slash K^{(l)}}
     \newcommand{\mh}{\mathfrak{M}_{\chi,\eta}}
     \newcommand{\mhaff}{\Phi_{K}^{-1}\paren{\eta}\slash K}
     \newcommand{\ycov}{\mathscr{Y}_{\chi,\eta}}
     \newcommand{\ot}{\cO\paren{T}}

     \newcommand{\ieta}{I_{\eta}}
     \newcommand{\ietal}{I^{\left(l\right)}_{\eta}}
     \newcommand{\el}{(l)}

    \newcommand{\brac}[1]{\ensuremath{\left\{#1\right\}}}
    \newcommand{\paren}[1]{\ensuremath{\left(#1\right)}}
    \newcommand{\sqbrac}[1]{\ensuremath{\left[#1\right]}}

    \newcommand{\gen}[1]{\ensuremath{\langle #1 \rangle}}

    \newcommand{\git}{\mathbin{
  \mathchoice{/\mkern-6mu/}
    {/\mkern-6mu/}
    {/\mkern-5mu/}
    {/\mkern-5mu/}}}
    \renewcommand\bar\overline
    \renewcommand\tilde\widetilde

    \newcommand{\eqn}[1]{\begin{eqnarray*}#1\end{eqnarray*}}
    \newcommand{\eqnn}[1]{\begin{eqnarray}#1\end{eqnarray}} 

    \newcommand{\theo}[1]{\begin{theorem}#1\end{theorem}}

    \newcommand{\prop}[1]{\begin{proposition}#1\end{proposition}}

    \newcommand{\lem}[1]{\begin{lemma}#1\end{lemma}}

    \newcommand{\cor}[1]{\begin{corollary}#1\end{corollary}}

    \providecommand{\defn}[1]{\begin{definition}#1\end{definition}}

    \newcommand{\rem}[1]{\begin{remark}#1\end{remark}}

    \renewcommand{\proof}[1]{\begin{myproof} #1\end{myproof}}

\usepackage{hyperref}
\hypersetup{
    colorlinks=true,       
    linkcolor=blue,          
    citecolor=blue,        
    filecolor=blue,      
    urlcolor=blue           
}

\begin{document}

\centerline{\bf \Large Quantum Multiplicative Hypertoric Varieties and Localization}
\vspace{0.5pc}\centerline{\sc Nicholas Cooney}

\parskip = 0pt

\begin{abstract} We consider $q$-deformations of multiplicative hypertoric varieties, where $q\in\bK^{\times}$ for $\bK$ an algebraically closed field of characteristic $0$. We construct an algebra $\qdiff$ of $q$-difference operators as a Heisenberg double in a braided monoidal category. We then focus on the case where $q$ is specialized to a root of unity. In this setting, we use $\qdiff$ to construct an Azumaya algebra on an $l$-twist of the multiplicative hypertoric variety, before showing that this algebra splits over the fibers of both the moment and resolution maps. Finally, we sketch a derived localization theorem for these Azumaya algebras.  \end{abstract}

\tableofcontents

\section{Introduction} \label{sec:intro}
This paper is concerned with deformations of multiplicative hypertoric varieties and, in particular, with the behaviour of these deformations when the deforming parameter is specialized to a root of unity. We will summarize our main results while simultaneously setting the scene, mentioning the literature which provided a point of departure for this work. Details of definitions that are contained in the main body of the text are omitted here.
\bigbreak
Hypertoric varieties were first introduced in \cite{BD} as Hyperk\"ahler quotients of smooth manifolds with torus actions. They can be described as GIT quotients of the cotangent bundle of affine $n$-space $T^{*}\bA^{n}$, with respect to a moment map $\mu:T^{*}\bA^{n}\rightarrow\mathfrak{t}^{*}$ induced by the action of the torus $T$ on $T^{*}\bA^{n}$, and with codomain the dual of the Lie algebra of the torus. Multiplicative hypertoric varieties, with which we are concerned here, can be defined analogously, but with a moment map $\Phi_{T}:T^{*}\bA^{n,\circ}\rightarrow T$ defined over an open subset of the cotangent bundle, and with codomain the group itself. Such maps were first defined and studied in \cite{AMM}. 
\bigbreak
In the second section, we discuss some preliminary results and definitions around the classical additive hypertoric varieties and their multiplicative counterparts. We also introduce the hypertoric enveloping algebra, first defined in \cite{MVdB} and more recently studied in \cite{BLPW}. We later define sheaves of differential operators on the multiplicative hypertoric variety which give a derived localization of central reductions of hypertoric enveloping algebras, in the sense of Beilinson and Bernstein's work for Lie algebras. 
\bigbreak
The third section is concerned with describing the key constructions required in later chapters, in the context of a braided monoidal category. After recalling the salient features of braided monoidal categories, and the notion of Hopf algebra objects in them, we give a general definition of a braided Heisenberg double in this context. We then go on to define a $q$-deformation of the Weyl algebra $\sD_{\bA^{n}}$ of differential operators on affine $n$-space. We construct such a deformation as a Heisenberg double of a deformation of the symmetric algebra and its dual, in the tensor category $\operatorname{Rep}_{q}\paren{T}$ of $\cO\paren{T}$-comodules, but with a non-trivial braiding defined by the grading induced from the $T$-action. We also introduce a deformation of the quantum multiplicative moment map, $\Phi_{q,T}:\cO\paren{T}\rightarrow\sD_{q}$. The definition of this deformed map is taken directly from the work of \cite{Jo}, in which the author defines a deformation quantization of multiplicative quiver varieties. Such quiver varieties where one takes the dimension at each vertex of the quiver to be $1$ are examples of the multiplicative hypertoric varieties that we focus on. 
\bigbreak
We finish by briefly sketching the definition of a category $\qdiff\paren{\mathfrak{M}_{q}}$ of \textit{quantum} $\sD$-\textit{modules on the quantum hypertoric variety} which, it is hoped, will lead to a treatment of $q$-deformations of multiplicative hypertoric varieties which unifies both the generic and root of unity cases. This attempt to provide a uniform definition of $\qdiff\paren{\mathfrak{M}_{q}}$ for all $q\in\bK^{\times}$ is informed by the papers \cite{BK1} and \cite{BK2}, in which the authors give a definition of both a quantum flag variety, deforming the category of quasi-coherent sheaves on the flag variety, and of quantum $\sD$-modules on the quantum flag variety, deforming the associated $\sD$-module category. The definition proceeds by categorically mimicking the steps involved in forming the GIT quotient that gives the multiplicative hypertoric variety, utilising the formalism of non-commutative projective schemes in \cite{AZ}.
\bigbreak
In section four, and for the rest of the document, we specialize to the case where the parameter $q$ is an $l^{th}$ root of unity. In this context, the behaviour of our deformed objects is radically different. This stems from the fact that the algebra of $q$-differential operators $\sD_{q}$ acquires a large centre, over which it is a finite dimensional module, and a sheaf over an $l$-twist of the cotangent bundle. We show that in this case $\sD_{q}$ is generically an Azumaya algebra over its centre. Recall that an Azumaya algebra on a variety is one which is locally a matrix algebra. We show that this is the case for $\sD_{q}$. This is a counterpart of the observation in \cite{BMR} that the sheaf of differential operators over affine $n$-space is an Azumaya algebra over its centre where the ground field has positive characteristic. 
\bigbreak
Indeed, one of the main motivations for this work is the analogy between geometry and representation theory over fields of positive characteristic, and over $q$-deformations for $q$ a root of unity. This is most evident in the analogy between results in \cite{BMR} and \cite{BK2} concerning, respectively, Lie algebras and $\sD$-modules on the flag variety, and Quantum Groups and Quantum $\sD$-modules on the quantum flag variety. In each case, there is a bounded derived equivalence between module categories for the relevant global sections algebras, and derived categories of $\sD$-modules. 
\bigbreak
In section five, we define a family of Azumaya algebras on an $l$-twist of the multiplicative hypertoric variety. The inspiration for this chapter is the paper \cite{BFG}, in which a general formalism is developed for quantum Hamiltonian reduction over fields of positive characteristic, and for constructing Azumaya algebras in this setting. In keeping with our analogy, we adapt their construction to define an Azumaya algebra on an $l$-twist of the multiplicative hypertoric variety $\ml$. We then show that this Azumaya algebra splits on fibers of the group-valued moment map $\ml\rightarrow H^{\el}$ and the map for which a multiplicative hypertoric variety is a symplectic resolution of singularities $\ml\rightarrow\Phi_{K^{\el}}^{-1}\paren{\eta^{l}}\slash K^{\el}$. A split Azumaya algebra is one which is globally a matrix algebra and so, by standard Morita theory, this induces an equivalence of abelian categories of $\sO$-modules on the fibers, and modules for the restriction of our Azumaya algebra to those fibers.
\bigbreak
In the final section, we sketch a derived localization theorem in the spirit of Beilinson and Bernstein obtaining, in the case where $\ml$ is smooth, a bounded derived equivalence between sheaves of modules over the Azumaya algebra and modules for a central reduction of a $q$-deformed multiplicative version of the hypertoric enveloping algebra. It is hoped that we can extend the Azumaya splitting of the previous chapter to formal neighbourhoods of the fibers, allowing for an algebro-geometric description of the representation theory of $q$-deformed hypertoric enveloping algebras. 

\paragraph{Acknowledgements.} This paper is an edited version of the author's doctoral thesis. The author would like to express his deepest gratitude to his doctoral supervisor, Kobi Kremnizer, for his generosity of spirit, mathematical insight, and friendship and, more specifically, for suggesting this project and for his guidance through it. This project was embarked upon independently by Iordan Ganev and he obtained the results presented here independently, in his article \cite{Ga}. The presentation of the material herein has benefited greatly from conversations with Ganev, and with David Jordan. This article can be read as something of a complement to Ganev's article, with the principal differences being the construction of $\sD_{q}$ as a braided Heisenberg double in section~\ref{sec:qdiff} and an alternative proof of theorem~\ref{adiffaz}. This work was supported by an EPSRC studentship and the EPSRC grant EP/I033343/1 in Motivic Invariants and Categorification.

\section{Preliminaries} \label{sec:prelim}

\subsection{Hypertoric Varieties}
In this work, we are concerned with quantizing the multiplicative versions of hypertoric varieties, which are quasi-Hamiltonian spaces arising from group-valued moment maps, first considered in \cite{AMM}. Let us first recall the more familiar additive hypertoric varieties, which were introduced in \cite{BD} as Hyperk\"ahler quotients of smooth manifolds with torus actions. 
\bigbreak
Henceforth, let $\bK$ be an algebraically closed field of characteristic 0, $T\cong\mathbb{G}^{n}_{m}$ be an $n$-torus for $n\in\mathbb{Z}_{>0}$, and $K\subset T$ a connected subtorus of rank $d\leq n$. Let $\phi:K\hookrightarrow T$ be an embedding given by $\phi\paren{k}_{i}=\prod_{j=1}^{d}k^{a_{ij}}$ for $a_{ij}\in\bZ$, $1\leq i\leq n$. The associated transpose map $\phi^{t}:T\rightarrow K$ is given by $\phi^{t}\paren{t}_{j}=\prod_{i=1}^{n}t^{a_{ji}}$. Let the algebras of functions on $T$ and $K$ be given by $\cO\paren{T}:=\bK\sqbrac{z_{1}^{\pm},\ldots,z_{n}^{\pm}}$ and $\cO\paren{K}:=\bK\sqbrac{u_{1}^{\pm},\ldots,u_{d}^{\pm}}$.
\bigbreak
Consider $T^{*}\mathbb{A}^{n}$ with the usual symplectic structure. The $T$-action on $\mathbb{A}^{n}$ differentiates to a map $\ft\rightarrow\text{Vect}\paren{\bA^{n}}$, where $\text{Vect}\paren{X}$ denotes the algebra of vector fields on an affine algebraic variety $X$. We also obtain an induced map $\fk\rightarrow\text{Vect}\paren{\bA^{n}}$, where $\mathfrak{k}:=\text{Lie}(K)$. This map extends uniquely to a map $\mu_{\mathscr{D}}:\text{Sym}\paren{\fk}\cong U\paren{\fk}\rightarrow\mathscr{D}\paren{\bA^{n}}$. This map, the \textit{quantum moment map} for the action of $K$ on $\bA^{n}$, is $K$-equivariant and filtered with respect to the canonical filtrations on $\text{Sym}\paren{\fk}$ and $\mathscr{D}\paren{\bA^{n}}$. The associated graded of this map, $\mu:=\text{gr}\paren{\mu_{\mathscr{D}}}:\text{Sym}\paren{\fk}\rightarrow\mathscr{O}_{T^{*}\bA^{n}}$, dualizes to the \textit{(additive) moment map} $$\mu: T^{*}\mathbb{A}^{n}\longrightarrow\mathfrak{t}^{*}\longrightarrow\mathfrak{k}^{*}$$ which factors through the moment map for the $T$-action. Let $\chi\in X^{*}(K)$ be a character of $K$ and $\lambda\in\mathfrak{k}^{*}$. We have the following:

\defn{The (additive) hyptertoric variety $\mathcal{M}_{\chi, \lambda}$ associated to the triple $\paren{K, \chi, \lambda}$ is the GIT quotient $$\mu^{-1}\paren{\lambda}\git_{\chi}K:=\mu^{-1}\paren{\lambda}^{\chi-ss}\git K$$ where $\mu^{-1}\paren{\lambda}^{\chi-ss}$ is the $\chi$-semistable locus of $\mu^{-1}\paren{\lambda}$ with respect to the $K$-equivariant line bundle $\mathscr{O}_{\mu^{-1}\paren{\lambda}}\paren{\chi}$.}

One can show that $\mathcal{M}_{\chi, \lambda}\cong\text{Proj}\paren{\bigoplus_{n\geq 0}\mathscr{O}\paren{\mu^{-1}\paren{\lambda}}^{\chi^{n}}}$ where $\mathscr{F}^{\chi^{n}}$ denotes the $\chi^{n}$-semi-invariant sections of a $K$-equivariant line bundle $\mathscr{F}$. 
\bigbreak
\subsection{Torus-valued moment maps}
In the paper \cite{AMM}, the authors define the notion of a \textit{quasi-Hamiltonian G-space} for $G$ a compact Lie group. This is a $G$-manifold $X$ where the moment map for the action takes values in the group $G$ itself, as opposed to the dual of $\fg:=\operatorname{Lie}(G)$. It is this notion of a group-valued moment map that we require for $G=T$, and which we later quantize. \\
Let $\paren{X,\omega}$ be a smooth symplectic variety equipped with a $T$-action preserving the symplectic form. Let $\gen{,}$ be a positive definite inner product by which we identify $\ft\cong\ft^{*}$ and, for $\xi\in\ft$, let $v_{\xi}$ be the corresponding vector field on $X$. Let $\theta\in\Omega^{1}(T,\ft)$ be the left-invariant Maurer-Cartan form. Since $T$ is abelian, this is equal to the right-invariant Maurer-Cartan form, which greatly simplifies the definition of the moment map.
\defn{A torus-valued moment map is a $T$-equivariant map $\Phi:X\rightarrow T$ such that $$\omega\paren{v_{\xi},-}=\gen{\Phi^{*}\theta,\xi}$$ for all $\xi\in\ft$.}
The corresponding map $\Phi^{\#}:\cO\paren{T}\rightarrow\cO\paren{X}$ is called the comoment map. 
The Maurer-Cartan form can be written locally as $\theta=\sum_{i} t^{-1}_{i}d_{i}$ which gives the moment map condition locally as $$\omega\paren{v_{\xi},-}=\sum_{i}\frac{d\Phi_{i}}{\Phi_{i}}\xi_{i}$$
The following two lemmas will be useful throughout this work. 
\lem{\label{mom-comp} Let $\phi:K\hookrightarrow T$ be an embedding of tori and let $\phi^{t}:T\rightarrow K$ denote its transpose. Suppose $\Phi:X\rightarrow T$ is a torus-valued moment map. Then $\phi^{t}\circ\Phi$ is a torus-valued moment map for the $K$-action induced by $\phi$.}
\proof{Let $\theta_{T}$ and $\theta_{K}$ denote the Maurer-Cartan forms of $T$ and $K$ respectively. Let $\operatorname{Lie}(\phi):\fk\rightarrow\ft$ and $\operatorname{Lie}(\phi^{t}):\ft\rightarrow\fk$ denote the linear maps on Lie algebras induced by $\phi:K\rightarrow T$ and $\phi^{t}:T\rightarrow K$. These are also transpose to one another. We use the same notation to denote the associated maps of $1$-forms $\operatorname{Lie}(\phi):\Omega^{1}\paren{-,\fk}\rightarrow\Omega^{1}\paren{-,\ft}$ and $\operatorname{Lie}(\phi^{t}):\Omega^{1}\paren{-,\ft}\rightarrow\Omega^{1}\paren{-,\fk}$. Also, for $\zeta\in\fk$, let $v^{K}_{\zeta}$ denote the corresponding vector field on $X$. The result follows from the computation:
\begin{equation*}
\begin{aligned}
\gen{\paren{\phi^{t}\circ\Phi}^{*}\theta_{K},\zeta}&=\gen{\Phi^{*}\circ\paren{\phi^{t}}^{*}\theta_{K},\zeta}=\gen{\Phi^{*}\circ\paren{\operatorname{Lie}(\phi^{t})\theta_{T}},\zeta}\\
&=\gen{\operatorname{Lie}(\phi^{t})\paren{\Phi^{*}\theta_{T}},\zeta}=\gen{\paren{\Phi^{*}\theta_{T}},\operatorname{Lie}(\phi)(\zeta)} \\
&=\omega\paren{v^{T}_{\operatorname{Lie}(\phi)(\zeta)},-}=\omega\paren{v^{K}_{\zeta},-}
\end{aligned}
\end{equation*}} 
\lem{\label{mom-quot} Suppose $\Phi:X\rightarrow T$ is a torus-valued moment map and that $K\subset T$ acts trivially on $X$. Then there is an induced action of $H:=T\slash K$ on $X$ and a moment map $\Phi_{H}:X\rightarrow H$.}
\proof{Let $\fh:=\operatorname{Lie}(H)$ and let $\psi:T\rightarrow H$ be the quotient map. There are short exact sequences of Lie algebras
\eqn{0\longrightarrow\fh\overset{\operatorname{Lie}(\psi^{t})}{\longrightarrow}\ft\overset{\operatorname{Lie}(\phi^{t})}{\longrightarrow}\fk\longrightarrow 0} and groups
\eqn{1\longrightarrow H\overset{\psi^{t}}{\longrightarrow}T\overset{\phi^{t}}{\longrightarrow}K\longrightarrow 1}
Given that $K$ acts trivially on $X$, the map $\phi^{t}\circ\Phi:X\rightarrow K$ is constant. Fix some $x^{\prime}\in X$ and let $t^{\prime}=\Phi(x^{\prime})$. We have that $\Phi(x)(t^{\prime})^{-1}\in\ker{\phi^{t}}$ for all $x\in X$. One can then define a map $\Phi_{H}:X\rightarrow H$ by $x\mapsto(\psi^{t})^{-1}(\Phi(x)(t^{\prime})^{-1})$. Let $L_{t}:T\rightarrow T$ denote left multiplication by $t\in T$. That $\Phi_{H}$ is a moment map can be verified using a similar computation to the previous lemma, and by noting that $\operatorname{Lie}(\phi)$ is surjective.
\begin{equation*}
\begin{aligned}
\omega\paren{v^{H}_{\operatorname{Lie}(\phi)(\xi)},-} &=\omega\paren{v^{T}_{\xi},-} 
=\gen{\Phi^{*}\theta_{T},\xi}=\gen{(\Phi^{*}_{H}\circ(\psi^{t})^{*}\circ L^{*}_{t^{\prime}}\theta_{T}, \xi}\\
&=\gen{\Phi^{*}_{H}(\operatorname{Lie}(\psi^{t})(\theta_{H})),\xi}=\gen{\operatorname{Lie}(\psi^{t})(\Phi^{*}_{H}\theta_{H}),\xi}=\gen{\Phi^{*}_{H}\theta_{H},\operatorname{Lie}(\psi)(\xi)}
\end{aligned}
\end{equation*}}
\subsection{The Multiplicative Hypertoric Variety}


We are now in a position to define the multiplicative hypertoric variety. Define an open subset $T^{*}\bA^{n,\circ}\subset T^{*}\bA^{n}$ by $T^{*}\bA^{n, \circ}:=\brac{\paren{p,w}| 1+p_{i}w_{i}\neq 0}\subset T^{*}\bA^{n}$. This is a symplectic variety equipped with a $T$-action induced by the usual scaling action. Locally, one can give a symplectic form on $T^{*}\bA^{n,\circ}$ by $$\omega=\sum_{i}\frac{dp_{i}\wedge dw_{i}}{\paren{1+p_{i}w_{i}}}$$ We have the following:
\prop{\label{tor-mom-map} The map $\Phi_{T}:T^{*}\bA^{n,\circ}\rightarrow T$ given by $\Phi_{T}:\paren{p,w}\mapsto\paren{1+p_{i}w_{i}}_{i}$ is a torus-valued moment map.}
\proof{We prove the result for $n=1$, proceeding by direct calculation. The general result follows similarly. Let $\xi\in\ft$.
\begin{equation*}
\begin{aligned}
\omega\paren{v_{\xi},-}=\frac{dp\wedge dw}{\paren{1+pw}}\paren{\xi p\frac{\partial}{\partial p}-\xi w\frac{\partial}{\partial w},-}=\xi\frac{wdp+pdw}{1+pw}=\xi\frac{d\Phi_{T}}{\Phi_{T}}=\gen{\Phi^{*}_{T}\theta_{T},\xi}
\end{aligned}
\end{equation*}}
Moreover, from Lemma~\ref{mom-comp} we have 
\cor{\label{tor-mom-map-K} The map $\Phi_{K}:T^{*}\bA^{n,\circ}\rightarrow K$ given by $$\Phi_{K}:\paren{p,w}\mapsto\paren{\prod^{n}_{i=1}\paren{1+p_{i}w_{i}}^{a_{ij}}}_{j}$$ is a torus-valued moment map.}
\defn{Given a triple $\paren{K,\chi,\eta}$ where $K\subset T$, $\chi\in X^{*}(K)$ and $\eta\in K$, the multiplicative hypertoric variety is defined as the GIT quotient $$\mh:=\Phi_{K}^{-1}(\eta)\git_{\chi}K$$}
As with the additive hypertoric variety, one can show that $$\mh\cong\operatorname{Proj}\paren{\bigoplus_{n\geq 0}\sO\paren{\Phi^{-1}_{K}(\eta)}^{\chi^{n}}}$$ 
We have the basic result
\prop{\label{mh-H-mom} Suppose $\mh$ is smooth, then it is a symplectic variety equipped with a torus-valued moment map $\Phi_{H}:\mh\rightarrow H$.}
\proof{There is an open affine cover of $\mh$ consisting of sets of the form $$U_{s}:=\operatorname{Spec}\paren{\sO(\Phi^{-1}_{K}(\eta))\sqbrac{s^{-1}}^{K}}$$ where $s\in\sO\paren{\Phi^{-1}_{K}(\eta)}^{\chi^{n}}$ for some $n\geq 0$. Each of these $U_{s}$ has a symplectic form with Hamiltonian $T$-action, where $K$ acts trivially. Hence, by Lemma~\ref{mom-quot}, there is a moment map $U_{s}\rightarrow H$ for each $U_{s}$. Where $U_{s}$ and $U_{r}$ have non-empty intersection, this is given by $U_{r}\cap U_{s}=U_{rs}$ and the inclusion of $U_{rs}$ into each of $U_{r}$ and $U_{s}$ is $H$-equivariant and the restriction of the moment maps agrees on this intersection. Hence, these glue together to give a moment map $\Phi_{H}:\mh\rightarrow H$.}

\subsection{\label{section:HEA} Quantum hamiltonian reduction and hypertoric enveloping algebras}
\label{subsection:HEA}
We now briefly discuss the connection between \textit{hypertoric Enveloping algebras}, first considered in \cite{MVdB}, and the non-commutative algebras we later construct, which quantize the coordinate rings of certain multiplicative hypertoric varieties, hence summarizing the content of the next few sections. These algebras are constructed via quantum Hamiltonian reduction of a $q$-deformation of the ring of polynomial differential operators on $\mathbb{A}^{n}$, with respect to a $K$-action and at a point $\eta\in K$. \bigbreak 
First, we define quantum Hamiltonian reduction more generally. Let $H$ be a Hopf $\bK$-algebra and $A$ an associative $\bK$-algebra in $H$-mod. A $\bK$-algebra homomorphism $\Phi:H\rightarrow A$ is called a \textit{quantum moment map} for the $H$-action on $A$, denoted by $\rhd:H\otimes A\rightarrow A$, if the following holds $$h\rhd a=\sum_{\paren{h}}\Phi\paren{h_{\paren{1}}}\cdot a\cdot\Phi\paren{S\paren{h_{\paren{2}}}}$$ for all $h\in H$ and $a\in A$. 

\rem{In the above definition, we have used \textit{Sweedler notation}, a way of simplifying expressions involving the coproduct of the Hopf algebra. Let $\Delta:H\rightarrow H\otimes H$ be the coproduct and let $h\in H$. There exist elements $h^{\paren{i}}_{\paren{1}},h^{\paren{i}}_{\paren{2}}\in H$ such that $$\Delta\paren{h}=\sum_{i}h^{\paren{i}}_{\paren{1}}\otimes h^{\paren{i}}_{\paren{2}}$$ In Sweedler notation, this is abbreviated to $\Delta\paren{h}=\sum_{\paren{h}}h_{\paren{1}}\otimes h_{\paren{2}}$}

Now let $I\trianglelefteq H$ be a Hopf ideal. As such, the quotient $A\slash A\cdot\Phi\paren{I}$ is equipped with an $H$-module structure. We have the following standard result:
\lem{The multiplication on $A$ descends to an algebra multiplication on the space of $H$-invariants $\paren{A\slash A\cdot\Phi\paren{I}}^{H}$.}
This algebra is called the \textit{quantum Hamiltonian reduction} of $A$ with respect to $H$ along $I$.
\bigbreak
The algebras for which we will later discuss localization are related to central reductions of \textit{hypertoric enveloping algebras}. The ring of polynomial differential operators on affine $n$-space, $\mathscr{D}\paren{\bA^{n}}$,  is isomorphic to the $n^{th}$--Weyl algebra $A_{n}$, the $\bK$-algebra generated by the elements $\brac{x_{1},\ldots,x_{n},\partial_{1},\ldots,\partial_{n}}$ with relations $\sqbrac{x_{i},x_{j}}=\sqbrac{\partial_{i},\partial_{j}}=0$ and $\sqbrac{\partial_{i},x_{j}}=\delta_{ij}$. This algebra is equipped with a $T$-action induced by the $T$-action on $\bA^{n}$ and is thus graded by the character lattice $X^{*}\paren{T}$ of $T$. 
We have the following algebra, introduced and studied in \cite{MVdB}.
\defn{The hypertoric enveloping algebra associated to $K\subset T$ is the algebra $U:=\mathscr{D}^{K}\paren{\bA^{n}}$ of $K$-invariant polynomial differential operators on $\bA^{n}$.}
The hypertoric enveloping algebra is equipped with a $T\slash K$-action and is hence graded by its character lattice $X^{*}\paren{T\slash K}$. It quantizes the categorical quotient $T^{*}\bA^{n}\git K$ and its central quotients quantize additive hypertoric varieties. The algebras $U$ and their quotients have a rich representation theory, as demonstrated in \cite{BLPW}, where a category $\cO$, having many properties analogous to the BGG category $\cO$, is considered. 
\bigbreak
As we have seen, the definition of multiplicative hypertoric varieties requires the passage to a certain open subset of $T^{*}\bA^{n}$ and, as such, any algebra quantizing this would involve localizing $\sD\paren{\bA^{n}}$. Unfortunately, the set of operators one would need to invert, $\brac{\paren{1+x_{i}\partial_{i}}^{l}|1\leq i\leq n, l\in\bZ_{>0}}$, does not form an Ore set in $\sD\paren{\bA^{n}}$. In order to remedy this, one considers a $q$-deformation $\sD_{q}\paren{\bA^{n}}$ of $\sD\paren{\bA^{n}}$ where $q\in\bK^{\times}$, in which the above set does satisfy the Ore condition. The algebra $\sD_{q}\paren{\bA^{n}}$ is generated by $\brac{x_{1},\ldots,x_{n},\partial_{1},\ldots,\partial_{n}}$ with relations $x_{i}x_{j}-qx_{j}x_{i}=\partial_{i}\partial_{j}-q\partial_{j}\partial_{i}=0$ for $i<j$ and $\partial_{i}x_{j}-qx_{j}\partial_{i}=(q-1)\delta_{ij}$. As we will observe in Section 3.2, the set of operators we wish to invert is Ore in $\sD_{q}\paren{\bA^{n}}$ and hence we can form the  localization $\mathscr{D}_{q}^{\circ}\paren{\bA^{n}}$. This too is graded by $X^{*}\paren{T}$ and also filtered by order of differential operators, with associated graded $\sO_{q}\paren{\paren{T^{*}\bA^{n}}^{\circ}}$, the corresponding localization of the $2n$-dimensional quantum plane. 
\bigbreak
As such, $q$-deformations of the algebras quantizing $\mathfrak{M}_{\chi,\eta}$ are quantum Hamiltonian reductions where we take $H=\mathcal{O}\paren{K}$, $A=\mathscr{D}^{\circ}_{q}\paren{\bA^{n}}$, $I=I_{\eta}$, the ideal generated by $\brac{x-\eta\cdot 1|x\in\cO\paren{K}}$, and $\Phi=\Phi_{q,K}$, a map related to the comoment map of the $K$-valued moment map from Proposition ~\ref{tor-mom-map}, which shall be introduced in Section 3.3. That is, $\cU_{q,\eta}:=\paren{\sD_{q}^{\circ}\paren{\bA^{n}}\slash\sD_{q}^{\circ}\paren{\bA^{n}}\cdot\Phi_{q,K}\paren{I_{\eta}}}^{K}$. They arise as central quotients of the $q$-deformed hypertoric enveloping algebra, $\sD_{q}^{K}\paren{\bA^{n}}$, localized to $T^{*}\bA^{n, \circ}$. We have the associated graded 
\eqn{
\begin{split}
\text{gr}\paren{\cU_{q,\eta}} &\cong\text{gr}\paren{\sD_{q}^{K}\paren{\bA^{n}}}^{\circ}\slash\text{gr}\paren{\sD_{q}^{K}\paren{\bA^{n}}}^{\circ}\cdot\Phi_{q,K} \paren{I_{\eta}} \\
&\cong\mathscr{O}_{q}\paren{\paren{T^{*}\bA^{n}}^{\circ}}^{K}\slash\gen{\paren{\Phi_{q,K}\paren{I_{\eta}}}} \\
&\cong\mathscr{O}_{q}\paren{\Phi^{-1}\paren{\eta}^{\circ}}^{K} 
\end{split}
}
We note that $\cU_{q,\eta}$ is hence finitely generated as an algebra over $\bK$, and that it has a Poisson structure inherited from the symplectic structure on $\mathfrak{M}_{\chi,\eta}$. As we shall see in the later sections of this work, where we impose the further restriction that $q$ is a primitive $l^{th}$ root of unity for $l>1$ odd, we can say much more: the algebra $\cU_{q,\eta}$ quantizes a twisted version of $\mh$, by Proposition ~\ref{locfree}; and we can demonstrate a derived equivalence between the categories of finitely generated modules for $\cU_{q,\eta}$ and of coherent sheaves of modules over an Azumaya algebra, constructed from $\sD_{q}\paren{\bA^{n}}^{\circ}$, over the twisted $\mh$.



\section{$\qdiff$ as a braided Heisenberg double} \label{sec:qdiff}
In this section, we construct the algebra $\qdiff$ of $q$-difference operators as a \textit{braided Heisenberg double}, a type of smash product algebra, in a braided monoidal category $\cC$. We begin by recalling the formalism of Hopf algebra objects in braided monoidal categories before defining the sheaf $\qdiff\paren{\bA^n}$ of $q$-deformed differential operators as a braided Heisenberg double, analogous to the classical construction of the Weyl algebra $\sD\paren{\bA^{n}}$ as a Heisenberg double, but in a category with non-trivial braiding. We finish by describing the $q$-deformed multiplicative quantum moment map in this setting, which is central to the constructions of the sheaf of algebras which quantizes $\mh$.  
\bigbreak
\subsection{Braided Hopf algebras and Heisenberg doubles}
We recall the definition of a Hopf algebra in a braided category, which closely follows the definition over the category of $\bK$-vector spaces. Let $\paren{\mathcal{C},\otimes,\sigma}$ be a braided monoidal category. Given two algebra objects $A$ and $B$ in $\cC$, one can define an algebra structure on $A\otimes B$ by defining the multiplication as
\begin{center}
\begin{tikzpicture}[node distance=7cm,on grid]
\node (A) at (0,0) {$m_{A,B}:A\otimes B\otimes A\otimes B$};
\node (B) [right=of A] {$A\otimes A\otimes B\otimes B$};
\node (C) [right=of B,xshift=-2cm] {$A\otimes B$};
\draw[->] (A) -- node [above] {$\text{id}\otimes\sigma_{B,A}\otimes\text{id}$} (B);
\draw[->] (B) -- node [above] {$m_{A}\otimes m_{B}$} (C);
\end{tikzpicture}
\end{center}
Hence, if we have a collection $A_{i}$, $1\leq i\leq n$, of algebra objects in $\cC$, we can inductively define a multiplication
\begin{center}
\begin{tikzpicture}[node distance=4cm,on grid]
\node (A) at (0,0) {$m_{\otimes_{i}A_{i}}:\bigotimes_{i}A_{i}\otimes\bigotimes_{i}A_{i}$};
\node (B) [right=of A] {$\bigotimes_{i}\paren{A_{i}\otimes A_{i}}$};
\node (C) [right=of B] {$\bigotimes_{i}A_{i}$};
\draw[->] (A) -- node [above] {} (B);
\draw[->] (B) -- node [above] {$\otimes m_{A_{i}}$} (C);
\end{tikzpicture}
\end{center}
The first map is given by a composition of the morphisms $\sigma$, $\sigma^{-1}$, the associator $\alpha:\paren{A_{i}\otimes A_{j}}\otimes A_{k}\rightarrow A_{i}\otimes\paren{A_{j}\otimes A_{k}}$ and its inverse $\alpha^{-1}$ (which we suppress in the notation) such that the order of the two copies of each factor $A_{i}$ is preserved. This composed morphism is unique by the coherence theorem of Joyal and Street \cite{JS}. We can now define
\defn{ A \textit{braided Hopf algebra} in $\cC$ is the data $\paren{H,m,\Delta,\eta,\epsilon,S}$ where $H$ is an algebra object in $\cC$ with multiplication $m:H\otimes H\rightarrow H$, and unit $\eta:\mathbbm{1}_{\cC}\rightarrow H$; a coalgebra with comultiplication $\Delta:H\rightarrow H\otimes H$ and counit $\epsilon:H\rightarrow\mathbbm{1}_{\cC}$; and a morphism $S:H\rightarrow H$, the antipode, such that the following diagram commutes:
\begin{center}
\begin{tikzpicture}[node distance=3cm,on grid]
\node (B1) at (0,0) {$H$};
\node (B2) [right=of B1] {$X$};
\node (B3) [right=of B2] {$H$};
\node (A1) [above=of B1,xshift=1.5cm,yshift=-1.33cm] {$H\otimes H$};
\node (A2) [right=of A1] {$H\otimes H$};
\node (C1) [below=of B1,xshift=1.5cm,yshift=1.33cm] {$H\otimes H$};
\node (C2) [right=of C1] {$H\otimes H$};
\draw[->] (B1) -- node [left] {$\Delta$} (A1);
\draw[->] (A1) -- node [above] {$\text{id}\otimes S$} (A2);
\draw[->] (A2) -- node [right] {$m$} (B3);
\draw[->] (B1) -- node [above] {$\epsilon$} (B2);
\draw[->] (B2) -- node [above] {$\eta$} (B3);
\draw[->] (B1) -- node [left] {$\Delta$} (C1);
\draw[->] (C1) -- node [below] {$S\otimes\text{id}$} (C2);
\draw[->] (C2) -- node [right] {$m$} (B3);
\end{tikzpicture}
\end{center}
} Analogous to the anti-homomorphic property of Hopf algebras, we have the following:
\lem{[{{\cite[Lemma 14.4]{Ma1}}}] The antipode morphism $S:H\rightarrow H$ of a braided Hopf algebra satisfies $$S\circ m=m\circ\sigma_{H,H}\circ\paren{S\otimes S}, \Delta\circ S=\paren{S\otimes S}\circ\sigma_{H,H}\circ\Delta$$
}
Dual to the definition of multiplication for a collection of algebra objects, one obtains a comultiplication on the tensor product of a collection of coalgebra objects $A_{i}$, $1\leq i\leq n$ in $\cC$
\begin{center}
\begin{tikzpicture}[node distance=4cm,on grid]
\node (A1) at (0,0) {$\Delta_{\otimes_{i}A_{i}}:\bigotimes_{i}A_{i}$};
\node (A2) [right=of A1] {$\bigotimes_{i}\paren{A_{i}\otimes A_{i}}$};
\node (A3) [right=of A2] {$\bigotimes_{i}A_{i}\otimes\bigotimes_{i}A_{i}$};
\draw[->] (A1) -- node [above] {$\otimes\Delta_{A_{i}}$} (A2);
\draw[->] (A2) -- node [above] {} (A3);
\end{tikzpicture}
\end{center}
where the second map is, as before, an unique morphism composed of braidings and associators. In a similar manner, one obtains a unit and counit of $\otimes_{i} A_{i}$ from those of the factors, and an antipode, using the properties of the previous lemma. We thus have 
\prop{Let $A_{i}$, $1\leq i\leq n$ be braided Hopf algebras in $\cC$. The tensor product $\bigotimes_{i}A_{i}$ is a braided Hopf algbera in $\cC$.}
\bigbreak
We make the following:
\defn{Let $A$ and $B$ be Hopf algebra objects in $\cC$. We say that there is a \textit{Hopf pairing} between $A$ and $B$ if there is an \textit{evaluation map} $\gen{-,-}:A\otimes B\rightarrow\mathbbm{1}_{\cC}$ compatible with the Hopf algebra structures as follows:
\begin{enumerate}
\item $\gen{-,-}\circ m_{A}=\paren{\gen{-,-}\otimes\gen{-,-}}\circ\paren{id\otimes\sigma_{A,B}\otimes id}\circ\Delta_{B}:A\otimes A\otimes B\rightarrow\mathbbm{1}_{\cC}$
\item $\gen{-,-}\circ m_{B}=\paren{\gen{-,-}\otimes\gen{-,-}}\circ\paren{id\otimes\sigma_{A,B}\otimes id}\circ\Delta_{A}:A\otimes B\otimes B\rightarrow\mathbbm{1}_{\cC}$
\item $\epsilon_{A}=\gen{-,-}\circ\eta_{B}:A\rightarrow\mathbbm{1}_{\cC}$,  $\epsilon_{B}=\gen{-,-}\circ\eta_{A}:B\rightarrow\mathbbm{1}_{\cC}$
\item $\gen{-,-}\circ S_{A}=\gen{-,-}\circ S_{B}:A\otimes B\rightarrow\mathbbm{1}_{\cC}$
\end{enumerate}
}
The existence of such a pairing exhibits $A$ as a left categorical dual of $B$. Note that, as opposed to the notion of a braided Hopf pairing given in \cite{Ma1}, we braid the middle two factors in 1.) and 2.), mimicking the convention in the symmetric monoidal case, where Majid first applies the pairing to the middle two factors, then the outer two. The absence of this transposition from Majid's definition results in a correspondence between braided right $B$-comodule algebras and braided left $A^{cop}$-module algebras\footnote{$A^{cop}$ is the Hopf algebra obtained from $A$ by setting $\Delta_{A^{cop}}:=\sigma^{-1}\Delta_{A}$.} in the category $\bar{\cC}$, that is, $\cC$ with the opposite braiding (see \cite{Lau}, Chapter 15 of \cite{Ma1} and Section 9.4 of \cite{Ma3}). 
\bigbreak
Now, suppose that $H$ and $H^{*}$ are dual Hopf algebras in $\cC$, with Hopf pairing $\gen{-,-}:H^{*}\times H\rightarrow\mathbbm{1}_{\cC}$ as above. The coproduct $\Delta:H\rightarrow H\otimes H$ equips $H$ with a right $H$-comodule structure. Mimicking the diagrammatic proof of Proposition 2.7 of \cite{Ta} for finite Hopf algebras we obtain a left action of $H^{*}$ on $H$, called the \textit{left regular action}, as follows:
\begin{center}
\begin{tikzpicture}[node distance=4cm,on grid]
\node (A1) at (0,0) {$act_{l}:H^{*}\otimes H$};
\node (A2) [right=of A1] {$H^{*}\otimes H\otimes H$};
\node (A3) [right=of A2] {$H^{*}\otimes H\otimes H$};
\node (A4) [right=of A3] {$H$};
\draw[->] (A1) -- node [above] {$\text{id}\otimes\Delta$} (A2);
\draw[->] (A2) -- node [above] {$\sigma_{23}$} (A3);
\draw[->] (A3) -- node [above] {$\gen{-,-}\otimes id$} (A4);
\end{tikzpicture}
\end{center}
The left regular action of $H^{*}$ on $H$, together with the action of $H$ on itself via left multiplication, leads to the following general definition:
\defn{\label{braided-Heis} Let $H$ be a dualizable Hopf algebra in $\cC$ and let $H'$ be a sub-Hopf algebra of its dual. The braided Heisenberg double $\mathfrak{h}\paren{H,H'}:= H\# H'$ is the smash product algebra with respect to the left regular action of $H'$ on $H$. That is, $\mathfrak{h}\paren{H,H'}$ is the object $H\otimes H'$ in $\cC$, equipped with multiplication 
\begin{center}
\begin{tikzpicture}[node distance=5cm,on grid]
\node (A1) at (0,0) {$H\otimes H'\otimes H\otimes H'$};
\node (A2) [right=of A1] {$H\otimes H'\otimes H'\otimes H\otimes H'$};
\node (A3) [right=of A2] {$H\otimes H'\otimes H\otimes H'\otimes H'$};
\node (B1) [below=of A1,yshift=4.2cm] {};
\node (B2) [right=of B1] {$H\otimes H\otimes H'\otimes H'$};
\node (B3) [right=of B2] {$H\otimes H'$};
\draw[->] (A1) -- node [above] {$\Delta_{2}$} (A2);
\draw[->] (A2) -- node [above] {$\sigma_{34}$} (A3);
\draw[->] (B1) -- node [above] {$act_{l,23}$} (B2);
\draw[->] (B2) -- node [above] {$m_{H}\otimes m_{H'}$} (B3);
\end{tikzpicture}
\end{center}
where $\Delta_{2}$ and $act_{l,23}$ indicates that these operators are applied to the second, and second and third, factors respectively, with the identity being applied elsewhere.}
Note that if the braiding is symmetric, one recovers the usual Heisenberg double multiplication. The algebra $\mathscr{D}_{q}\paren{\bA^{n}}$ will be described as such a braided Heisenberg double.

\subsection{$q$-Deformed Differential Operators}
Let $M=\sqbrac{m_{ij}}$ be an $n\times n$ matrix of indeterminates such that $m_{ji}=-m_{ij}$ for all $i\neq j$ and let $Q=\sqbrac{q_{ij}}$ be the matrix with entries $q_{ij}:=q^{m_{ij}}$ for some indeterminate $q$, which we later specialize to values in $\bK^{\times}$. Let $\text{Rep}_{q}\paren{T}$ be the category of $\ot$-comodules with the following braiding:

\eqn{\sigma_{V,W}:V\otimes W\rightarrow W\otimes V}
\eqnn{\label{braid} v\otimes w\mapsto q^{\text{deg}\paren{v}\cdot M\cdot\text{deg}\paren{w}^{t}}w\otimes v}

Note that, for $q=1$, we recover the category of $\cO\paren{T}$-comodules $\operatorname{Rep}\paren{T}$. Let $\Lambda$ be the cocharacter lattice of $T$ and let $\Lambda_{\bK}=\Lambda\otimes_{\bZ} \mathbb{K}$. Denote by $T\paren{\Lambda_{\bK}}$ its tensor algebra. Let $S_{q}\paren{\Lambda_{\bK}}$ be the algebra 

\eqn{T\paren{\Lambda_{\bK}}\slash\gen{v_{i}\otimes w_{j}-\sigma\paren{v_{i}\otimes w_{j}}, i\neq j}} where $v=\paren{v_{1},\ldots,v_{n}}\in \Lambda_{\bK}$.
\rem{$S_{q}\paren{\Lambda_{\bK}}$ is isomorphic to the quantum plane as an algebra object in $\text{Vect}_{\bK}$. It is not, however, a $\sigma$-commutative algebra object in $\text{Rep}_{q}\paren{T}$, as we do not impose the corresponding relation when $i=j$ and so do not take the $q$-symmetric algebra as our analogue of $S\paren{\Lambda_{\bK}}$ here. This would require $q_{ii}^{2}=1$ which, as we shall see when forming the Heisenberg double, would give the Heisenberg relation $\partial_{i}x_{i}-x_{i}\partial_{i}-1$ as opposed to the $q$-deformed version we require.}
\prop{Let $V=\bK x$ be a one-dimensional vector space. The braided tensor algebra $T\paren{V}=\bK\gen{x}$ is an algbera with coproduct, counit, and antipode given on generators by $\Delta:x\mapsto 1\otimes x + x\otimes 1$, $\epsilon\paren{x}=0$ and $S\paren{x}=-x$ respectively.}
\proof{Consider the diagonal embedding $V\rightarrow V\oplus V$. Apply the tensor algebra functor to obtain a map $T(V)\rightarrow T(V\oplus V)$. There is an isomorphism $T(V\oplus V)\cong T(V)\otimes T(V)$ given on generators by $\paren{x,y}\mapsto x\otimes 1 + 1\otimes y$. By the universal property of the tensor algebra, it is enough to check the Hopf algebra axioms for elements of $V$ which, by the above, are clear.}
\bigbreak
Let $\brac{x_{1},\ldots,x_{n}}$ be a basis for $\Lambda_{\bK}$ and hence a generating set for $S_{q}\paren{\Lambda_{\bK}}$. By examining the presentation given for $S_{q}\paren{\Lambda_{\bK}}$ above, one can see that $S_{q}\paren{\Lambda_{\bK}}\cong\bigotimes_{i=1}^{n}\bK\gen{x_{i}}$ as an algebra object in $\text{Rep}_{q}\paren{T}$. As such, $S_{q}\paren{\Lambda_{\bK}}$ is a braided tensor product of one-dimensional tensor algebras and so we have:
\cor{The algebra $S_{q}\paren{\Lambda_{\bK}}$ is a braided Hopf algebra in $\text{Rep}_{q}\paren{T}$.}
\bigbreak
Let $\Lambda^{*}_{\bK}=\Lambda^{*}\otimes_{\bZ}\bK$. In a similar manner, we form the algebra $S_{q}\paren{\Lambda^{*}_{\bK}}$. These algebras are dual braided Hopf algebras under the pairing induced by the pairing on dual vector spaces and, as such, we can form their Heisenberg double $\mathfrak{h}\paren{S_{q}\paren{\Lambda_{\bK}},S_{q}\paren{\Lambda^{*}_{\bK}}}$. The following algebra is obtained via direct computation from the multiplication defined in Definition~\ref{braided-Heis}:
\prop{Let $\brac{x_{1},\ldots,x_{n}}$ be a basis for $\Lambda_{\bK}$ and generating set for $S_{q}\paren{\Lambda_{\bK}}$ and let $\brac{\partial_{1},\ldots,\partial_{n}}$ be a basis for $\Lambda^{*}_{\bK}$, dual to that for $\Lambda_{\bK}$, and a generating set for $S_{q}\paren{\Lambda^{*}_{\bK}}$. The Heisenberg double $\mathfrak{h}\paren{S_{q}\paren{\Lambda_{\bK}},S_{q}\paren{\Lambda^{*}_{\bK}}}$ is the algebra isomorphic to $S_{q}\paren{\Lambda_{\bK}}\otimes S_{q}\paren{\Lambda^{*}_{\bK}}$ as a vector space, with relations
\begin{align*}
x_{i}\otimes x_{j} - q_{ij}x_{j}\otimes x_{i} &=0 \\
\partial_{i}\otimes\partial_{j} - q_{ij}\partial_{j}\otimes\partial_{i} &=0 \\
\end{align*} for all $i\neq j$ and 
\eqn{\partial_{i}\otimes x_{j}-q^{-1}_{ij}x_{j}\otimes\partial_{i}-\delta_{ij}=0} for all $1\leq i,j \leq n$.}

Where there is no risk of ambiguity, we denote $S_{q}\paren{\Lambda_{\bK}}$ by $S_{q}$ and $S_{q}\paren{\Lambda^{*}_{\bK}}$ by $S^{*}_{q}$. 
Alternatively, we can describe this algebra as a braided tensor product of copies of $A_{q}=\bK\gen{x,\partial}\slash\paren{\partial x-qx\partial=1}$, the first $q$-Weyl algebra. That is,
\eqn{\mathfrak{h}\paren{S_{q},S^{*}_{q}}\cong\bigotimes_{1\leq i\leq n}A_{q^{-1}_{ii}}} with the braiding between coordinates giving rise to the remaining relations. 
\bigbreak
Given a matrix $M$ as above, define a matrix $M^{\prime}$ with coefficients $m^{\prime}_{ij}=m_{ij}$ for $i\neq j$ and $m^{\prime}_{ii}=-m_{ii}$ for all $i$, labelling $q^{\prime}_{ij}:=q^{m^{\prime}_{ij}}$. In order that our algebra better resembles algebras of $q$-difference operators found elsewhere in the literature, we perform the above construction again, but with reference to $M^{\prime}$. We also rescale the algebra by replacing each $x_{i}$ with $x_{i}\slash\paren{q^{\prime}_{ii}-1}$. This produces an algebra isomorphic to $\mathfrak{h}\paren{S_{q},S^{*}_{q}}$ for $q\neq 1$ and gives the \emph{$q$-Heisenberg relation}
\eqn{\partial_{i}\otimes x_{j}-q^{\prime}_{ij}x_{j}\otimes\partial_{i}-\delta_{ij}\paren{q^{\prime}_{ii}-1}=0}This is the algebra we take as our $\qdiff\paren{\bA^{n}}$. The algebras $S_{q}$ and $\qdiff\paren{\bA^{n}}$ are flat deformations of the corresponding classical objects, with the same bases of monomials over the coefficient ring $\bK\sqbrac{\sqbrac{(q^{\prime}_{ij})_{1\leq i, j\leq n}}}$. 
Henceforth, for ease of notation, we again just write $q_{ij}$ in place of $q^{\prime}_{ij}$. Straightforward calculations give the following:
\lem{For $1\leq i\leq n$, let $\alpha_{i}:=1+x_{i}\partial_{i}$ be the \textit{Euler operator}. The braided algebra $\qdiff\paren{\bA^{n}}$ satisfies the following identities:
\begin{enumerate}
\item $\partial_{i} x_{i}^{n}=q_{ii}^{n}x_{i}^{n}\partial_{i}+\paren{q_{ii}^{n}-1}x_{i}^{n-1}$ and $\partial_{i}^{n}x_{i}=q_{ii}^{n}x_{i}\partial_{i}^{n}+\paren{q_{ii}^{n}-1}\partial_{i}^{n-1}$
\item $\alpha_{i}x_{j}=q^{\delta_{ij}}_{ii}x_{j}\alpha_{i}$ and $\alpha_{i}\partial_{j}=q^{-\delta_{ij}}_{ii}\partial_{j}\alpha_{i}$
\end{enumerate}
} The $q$-commutativity of the $\alpha_{i}$'s means that the set $\brac{\alpha_{i}| 1\leq i\leq n}$ satisfies the Ore condition and we can form the localization $\qdiff^{\circ}=\qdiff\sqbrac{(\prod\alpha_{i})^{-1}}$ which will appear again in the next section.
\subsection{The moment map}
Before coming to a description of the categories of interest, it is necessary to discuss what a moment map is in this braided context, similar to the definition in section 2.4. 
\newline
Let $H$ be a Hopf algebra in a braided category $\paren{\cC,\otimes,\sigma}$. We have the following definition, analogous to that in a symmetric monoidal category:
\defn{An algebra object $\paren{A,m_{A},\eta_{A}}$ in $\cC$ is called an $H$-module algebra if it is equipped with an $H$-module structure\footnote{The morphism $\rhd: H\otimes A\rightarrow A$ equips $A$ with an $H$-module structure if $$\rhd\circ \paren{m_{H}\otimes id}=\rhd\circ\paren{id\otimes\rhd}:H\otimes H\otimes A\rightarrow A$$} $\rhd: H\otimes A\rightarrow A$ satisfying the following conditions
\begin{enumerate}
\item $\rhd\circ\paren{id_{H}\otimes m_{A}}=m_{A}\circ\paren{\rhd\otimes \rhd}\circ\sigma_{23}\circ\Delta_{H}:H\otimes A\otimes A\rightarrow A$
\item $\rhd\circ\eta_{A}=\eta_{A}\circ\epsilon_{H}:H\rightarrow A$
\end{enumerate}
} Let $A$ be such an algebra. A homomorphism $\Phi:H\rightarrow A$ is called a quantum moment map for this action if, for all $h\in H$ and $a\in A$, the action is given by the composition of maps
\begin{center}
\begin{tikzpicture}[node distance=4.5cm,on grid]
\node (A1) at (0,0) {$H\otimes A$};
\node (A2) [right=of A1] {$H\otimes H\otimes A$};
\node (A3) [right=of A2] {$H\otimes A\otimes H$};
\node (A4) [right=of A3] {};
\node (B1) [below=of A1,yshift=3.7cm] {};
\node (B2) [right=of B1] {$H\otimes A\otimes H$};
\node (B3) [right=of B2] {$A\otimes A\otimes A$};
\node (B4) [right=of B3] {$A$};
\draw[->] (A1) -- node [above] {$\Delta\otimes\text{id}$} (A2);
\draw[->] (A2) -- node [above] {$\text{id}\otimes\sigma_{H,A}$} (A3);
\draw[->] (A3) -- node [above] {$\text{id}\otimes\text{id}\otimes S$} (A4);
\draw[->] (B2) -- node [above] {$\Phi\otimes\text{id}\otimes\Phi$} (B3);
\draw[->] (B3) -- node [above] {$m_{A}\circ\paren{\text{id}\otimes m_{A}}$} (B4);
\end{tikzpicture}
\end{center}
\bigbreak
Consider $\cO\paren{T}$ as an $\cO\paren{T}$-module under the adjoint action, so that all elements have degree $0$ and, as such, the adjoint action is trivial. By the relations in the previous lemma, the map $\Phi_{q,T}:\cO\paren{T}\rightarrow \mathscr{D}_{q}$ given by $\Phi_{q,T}:z_{i}\mapsto 1+x_{i}\partial_{i}$ is a moment map for the action of $\cO\paren{T}$ on $\mathscr{D}_{q}$, agreeing with the grading induced by the $\cO\paren{T}$-coaction. 
\bigbreak
The embedding $\phi:K\hookrightarrow T$ given by $\phi\paren{k}_{i}=\prod_{j=1}^{d}k_{j}^{a_{ij}}$ induces a map, the pull-back of the transpose of $\phi$, $\phi^{t,*}: \cO\paren{K}\rightarrow\cO\paren{T}$ given by $$u_{j}\mapsto\prod^{n}_{i=1}z_{i}^{a_{ij}}$$ where $\cO\paren{K}=\bK\sqbrac{u_{1}^{\pm},\ldots,u_{d}^{\pm}}$. The composition
\eqn{\Phi_{q,K}:\cO\paren{K}\overset{\phi^{t,*}}{\longrightarrow}\cO\paren{T}\overset{\Phi_{q,T}}{\longrightarrow}\qdiff} given by 
\eqn{u_{j}\mapsto\prod_{i=1}^{n}\paren{1+x_{i}\partial_{i}}^{a_{ij}}}
is a moment map for the corresponding $\cO\paren{K}$-coaction. The map $\Phi_{q,K}$ both deforms and quantizes the comoment map associated to the torus-valued moment map $\Phi_{K}:T^{*}\bA^{n,\circ}\rightarrow K$ introduced in Corollary 2.6.
\remark{Using the constructions of the previous sections, it should be possible to define a category $q$-deforming the category of $\mathscr{D}$-modules on the multiplicative hypertoric variety for an arbitrary $q\in\bK^{*}$, by categorically mimicking the steps required to construct $\mh$. That is, we can define a category $\qdiff\paren{\mathfrak{M}_{q}}$ of \textit{quantum} $\sD$-\textit{modules on the quantum multiplicative hypertoric variety} as follows. Let $X=T^{*}\bA^{n,\circ}$ with an action of $K\subset T$ as before. Let $\ieta\trianglelefteq\cO\paren{K}$ be the ideal generated by the elements $\brac{u_{j}-\eta_{j}}_{1\leq j\leq d}$ for $\eta\in K$. For $\chi\in X^{*}\paren{K}$, let $\sO\paren{\chi}$ be the associated $K$-equivariant line bundle. We define $\qdiff\paren{\mathfrak{M}_{q}}$ to be the quotient of the category $\cD$ of $K$-equivariant $\qdiff$-modules on $X$ annihilated by $\Phi_{q,K}\paren{\ieta}$, by the category of modules $M$ in the kernel of the functor $$\bar{\Gamma}:M\mapsto\bigotimes_{-\infty}^{\infty}Hom_{\cD}\paren{\sO_{X},M\otimes\sO\paren{\chi}^{\otimes n}}$$
(see \cite{AZ} where this formalism of noncommutative projective schemes is developed). It is expected that the category we have just sketched for an arbitrary $q\in\bK^{*}$ will, when $q$ is specialized to an $l^{th}$ root of unity, allow us to recover the category of modules for the Azumaya algebra over the $l$-twisted hypertoric variety defined and studied in section 5.}



\section{$\sD_{q}$ at a Root of Unity} \label{sec:ltwist}
We now proceed to the study of $\qdiff$ in the case where we specialize $q$ to an $l^{th}$ root of unity. Here, $\qdiff$ becomes a sheaf of algebras over an $l$-twist of the cotangent bundle $T^{*}\bA^{n}$. It attains a large centre, over which it is a finite dimensional module and, in fact, generically an Azumaya algebra. After recalling the definition and basic properties of Azumaya algebras, we summarize the basic structure of $\qdiff$ for $q$ a root of unity, with the key property being that it is an Azumaya algebra over an affine open subset of an $l$-twist of the cotangent bundle on affine $n$-space.
\subsection{Azumaya algebras}
The main reference for the results in this section is \cite{Mi}. Let $R$ be a commutative ring. An $R$-algebra $A$ which is finitely generated and projective as an $R$-module, is called \textit{Azumaya} if the map
\eqn{A\otimes_{R} A^{opp}\longrightarrow End_{R}\paren{A}}
is an isomorphism. Azumaya algebras can also be characterised as those which, locally, are central simple algebras. 
\prop{Let $A$ be an $R$-algebra as above. Then $A$ is an Azumaya algebra if, and only if, $A\otimes_{R}k\paren{p}$ is a central simple $k\paren{p}$-algebra for all closed points $p\in Spec\paren{R}$.}
Now, let $X$ be a scheme and $\mathscr{A}$ be a locally free sheaf of $\mathscr{O}_{X}$ algebras finite rank on $X$. The sheaf $\mathscr{A}$ is \textit{Azumaya} if the map
\eqn{\sA\otimes_{\sO_{X}}\sA^{opp}\longrightarrow\mathscr{E}nd_{\sO_{X}}\paren{\sA}}
is an isomorphism. If one now assumes that $X$ is connected and of finite type over an algebraically closed field, then $\sA$ is a sheaf of Azumaya algebras precisely when $\sA$ has constant rank on the set of closed points and $\sA\otimes_{\sO_{X}}k\paren{p}$ is a matrix algebra for any closed point $p\in X$. Finally, we say that a sheaf of Azumaya algebras is \textit{split} if it is globally a matrix algebra: that is, if there exists a vector bundle $\mathscr{V}$ on $X$ such that $\sA\cong\mathscr{E}nd_{\sO_{X}}\paren{\mathscr{V}}$ on $X$. If $\sA$ is a split Azumaya algebra on $X$, then $\sA$ is Morita equivalent to $\sO_{X}$. 
The following technical result will be used later when we discuss Azumaya splitting on hypertoric varieties:
\lem{\label{Cart-Az} Let $\sA$ be an Azumaya algebra on $Y$ and let 
\begin{center}
\begin{tikzpicture}[node distance=2cm,on grid]
\node (X) at (0,0) {$X$};
\node (Y) [right=of X] {$Y$};
\node (U) [below=of X] {$U$};
\node (V) [right=of U] {$V$};
\draw[->] (X) -- node [above] {$f$} (Y);
\draw[->] (X) -- node [left]  {$s$} (U);
\draw[->] (Y) -- node [right] {$t$} (V);
\draw[->] (U) -- node [below] {$g$} (V);
\end{tikzpicture}
\end{center}
be a Cartesian square such that $f:X\rightarrow Y$ is a covering map and $f^{*}\sA$ is split on $X$. Then $\sA$ splits when restricted to the fibers of $t:Y\rightarrow V$.}
\proof{Let $v\in V$ such that $t^{-1}(v)$ is non-empty and let $y\in t^{-1}(v)$. Since $f$ is a cover, there exists some $x\in X$ such that $f(x)=y$ and, hence, some $u\in U$ such that $g(u)=v$. Since the square is Cartesian, the restriction of $f$ gives an isomorphism of schemes $f^{\prime}:=f\vert_{s^{-1}(u)}:s^{-1}(u)\rightarrow t^{-1}(v)$. Since $f$ splits $\sA$, there is a vector bundle $\sV$ on $X$ such that $f^{*}\sA\cong\mathscr{E}nd_{\sO_{X}}\paren{\sV}$. Hence, we have that $\sA\vert_{t^{-1}(v)}\cong(f^{\prime})_{*}(f^{\prime})^{*}\sA\vert_{t^{-1}(v)}\cong\mathscr{E}nd_{\sO_{t^{-1}(v)}}\paren{(f^{\prime})_{*}\sV\vert_{s^{-1}(u)}}$.}

\subsection{$\qdiff$ as an Azumaya algebra}
Let $l>1$ be an odd integer and let $q$ be a primitive $l^{th}$ root of unity. 
For $X=\bA^{n}$, we define $T^{*}X^{\paren{l}}\cong Spec\paren{\bK\sqbrac{x^{l}_{1}\ldots,x^{l}_{n},\partial^{l}_{1}\ldots,\partial^{l}_{n}}}$, the $l$-twist of the cotangent bundle of $X$. The embedding $\bK\sqbrac{x^{l}_{1}\ldots,x^{l}_{n},\partial^{l}_{1}\ldots,\partial^{l}_{n}}\hookrightarrow\bK\sqbrac{x_{1}\ldots,x_{n},\partial_{1}\ldots,\partial_{n}}$ induces a finite map of $\bK$-schemes, which we denote $Fr_{l}:T^{*}X\rightarrow T^{*}X^{\paren{l}}$.
\bigbreak
Henceforth, for simplicity, we take a single parameter deformation of the Weyl algebra: that is, with previous notation, $m_{ij}=1$ for all $1\leq i,j\leq n$. For $q$ a primitive $l^{th}$ root of unity, $\qdiff$ contains as a central subalgebra its \textit{$l$-centre} $Z_{l}\paren{\qdiff}=\bK\sqbrac{x^{l}_{1}\ldots,x^{l}_{n},\partial^{l}_{1}\ldots,\partial^{l}_{n}}$. The algebra $\qdiff$ is a free module of rank $l^{2n}$ over $Z_{l}\paren{\qdiff}$, with a basis given by $$\brac{x^{r}_{i}\partial^{s}_{j}|1\leq i,j\leq n,0\leq r,s\leq l-1}$$ We have, moreover
\lem{\label{powerlem} The centre $Z\paren{\qdiff}$ of $\qdiff$ is given by $Z_{l}\paren{\qdiff}$.} 
\proof{We demonstrate the result for $n=1$, from which the general case straightforwardly follows. Clearly, $Z_{l}\paren{\qdiff}\subset Z\paren{\qdiff}$. Let $z=\sum_{m,n}c_{m,n}x^{m}\partial^{n}\in Z\paren{\qdiff}$. Given that $\partial z=z\partial$ and $xz=zx$, one can determine that $$c_{m,n}\paren{1-q^{m}}=\paren{-1}^{i}c_{m+i,n+i}\paren{1-q^{m+i}}$$ and $$c_{m,n}\paren{1-q^{n}}=\paren{-1}^{i}c_{m+i,n+i}\paren{1-q^{n+i}}$$ for $i\geq 0$. From this, one can see immediately that $c_{m,n}=0$ when $l\nmid m,n$.} \bigbreak Set $\Delta=\prod_{i=1}^{n}\alpha_{i}$. 
\lem{\label{delta-l}We have $\Delta^{l}=\prod_{i=1}^{n}1+x_{i}^{l}\partial_{i}^{l}$.}
\proof{Again, it suffices to demonstrate the result for $n=1$, where $\Delta=\alpha=1+x\partial$. Since $\alpha$ conjugates elements by a root of unity, $\alpha^{l}\in Z\paren{\qdiff}=Z_{l}\paren{\qdiff}$ and we have $\alpha^{l}=a+bx^{l}+c\partial^{l}+dx^{l}\partial^{l}$. One can readily see that $a=1$. Since $Z_{l}\paren{\qdiff}\subset k\sqbrac{x_i,\partial_i}$ is a graded subalgebra with $\operatorname{deg}\paren{\Delta}=0$, we have that $\operatorname{deg}\paren{\Delta^{l}}=0$ and so $b=c=0$. We claim that for $l$ odd, $d=q^{l(l-1)\slash 2}$. Proceeding by induction, the statement is true for $l=1$. Suppose it is true for some odd $l>1$. Then
\begin{align*}
q^{l(l-1)\slash 2}x^l\partial^l x\partial x\partial &=q^{l(l-1)\slash 2}x^l\paren{q^{l}x\partial^l+\paren{q^{l}-1}\partial^{l-1}}\partial x\partial \\
&=q^{l+l(l-1)\slash 2}x^{l+1}\partial^{l+1}x\partial+q^{l(l-1)\slash 2}\paren{q^{l}-1}x^{l}\partial^{l}x\partial \\
&=q^{l+l(l-1)\slash 2}x^{l+1}\paren{q^{l+1}x\partial^{l+1}+\paren{q^{l+1}-1}\partial^{l}}\partial+\ldots \\
&=q^{2l+1+l(l-1)\slash 2}x^{l+2}\partial^{l+2}+\ldots \\
&=q^{(l+2)(l+1)\slash 2}x^{l+2}\partial^{l+2}+\ldots \\
\end{align*} thus establishing the claim. For $l>1$ odd and $q$ a primitive $l^{th}$ root of unity, $q^{l(l-1)\slash 2}=1$ and so we are done.}
\bigbreak
We have the following (\cite{Ba}, Proposition 2.3)
\theo{\label{Azum-dq} The Azumaya locus of $\qdiff$ is $U_{\Delta^{l}}\subset Spec\paren{Z_{l}\paren{\qdiff}}$ where $$U_{\Delta^{l}}=\brac{\mathfrak{m}\in\operatorname{MaxSpec}\paren{{Z_{l}\paren{\qdiff}}}|\Delta^{l}\notin\mathfrak{m}}$$}
In order to prove this theorem, we first demonstrate it for $\qdiff$ where $n=1$, and then show that the result holds for a braided tensor product of finitely many copies of it. We sketch the proof for $n=1$ from \cite{Ba}, where the braiding is trivial.
\bigbreak
\proof{Let $n=1$ and let $\zeta=\paren{a,\omega}\in Spec\paren{Z_{l}}$. Suppose there is an isomorphism $\phi:\mathscr{D}_{q,\zeta}\overset{\sim}\rightarrow M_{l}\paren{\bK}$ and let $X:=\phi(x)$ and $Y:=\phi(\partial)$. Suppose first that $a\neq 0$. Since $X^{l}=a$ we have that $X=diag\paren{\lambda_{1},\ldots,\lambda_{l}}$ where $\lambda^{l}_{i}=a$. Since the elements $I,X,\ldots,X^{l-1}$ are linearly independent, we have $\lambda_{i}\neq\lambda_{j}$ for $i\neq j$. After permuting the basis elements, we have $\lambda_{i}=\lambda\cdot q^{i}$. The equation $XY-qYX=q-1$ gives 
 $$Y=
\begin{pmatrix}
  -1\slash\lambda_{0} & 0 & 0 & \cdots & b_{0,l-1} \\
  b_{1,0} & -1\slash\lambda_{1} & 0 & \cdots & 0 \\
  0  & b_{2,1} & -1\slash\lambda_{2} & \cdots & 0  \\
  \vdots & \vdots & \vdots & \ddots & \vdots \\
  0 & 0 & \cdots & b_{l-1,l-2} & -1\slash\lambda_{l-1}
 \end{pmatrix}$$ for some $b_{ij}\in\bK$. One then checks that $\bK^{l}$ is an irreducible $\mathscr{D}_{q,\zeta}$-module if and only if each $b_{ij}\neq 0$. Given that $Y^{l}=\omega$, each $b_{ij}\neq 0$ if, and only if, $1+a\omega\neq 0$. One proves the $\omega\neq 0$ case similarly. \newline It remains to prove the $a=\omega=0$ case. Here, we have that $\bK\sqbrac{x}\slash\paren{x^{l}}$ is an irreducible $\sD_{q,\zeta}$-module and we are done.}
Since matrix algebras are semi-simple, one would expect that any $q$-deformations of these algebras would be isomorphic to the non-deformed algebra. The following technical lemma is a version of that proved in Proposition 2.27 of \cite{Ga}.
\lem{Let $\paren{\cC,\otimes,\sigma}$ be a braided tensor category over a braided commutative $k$-algebra $R$, where the braiding is determined by the grading as in equation~\ref{braid}. A tensor product of matrix algebras in $\cC$ is a matrix algebra, isomorphic to the tensor product of matrix algebras with trivial braiding.} As such, to demonstrate that a given algebra is Azumaya in the braided setting, it suffices to show that it is Azumaya in the non-deformed setting, as shown above.
\bigbreak
Since $\Delta^{l}\in Z_{l}\paren{\qdiff}$, we can form the localization $Z_{l,\Delta^{l}}$ and obtain $$\qdiff^{\circ}:=\qdiff\otimes_{Z_{l}}Z_{l,\Delta^{l}}$$ Note that $Spec\paren{Z_{l,\Delta^{l}}}\cong T^{*}\bA^{n,(l),\circ}$. Finally, from the above calculations, we have
\cor{The algebra $\qdiff^{\circ}$ is Azumaya over its centre $Z_{l,\Delta^{l}}$.}

\section{An Azumaya algebra quantizing the Multiplicative Hypertoric Variety} \label{sec:azumaya}
In the paper \cite{BFG}, a formalism for quantum Hamiltonian reduction over fields of positive characteristic is developed and used to construct Azumaya algebras on a Frobenius twist of the Hilbert scheme of points in the plane. Here, we adopt this formalism in our setting to construct an Azumaya algebra $\adiff$ on an $l$-twist of the hypertoric variety, $\fM^{\paren{l}}_{\chi,\lambda}$, as a quantum Hamiltonian reduction of $\qdiff$. We then show that this algebra splits on the fibers of the resolution map $$\Psi:\ml\longrightarrow\Phi_{K^{\el}}^{-1}\paren{\eta^{l}}\slash K^{\el}$$ and of the moment map $$\Phi_{H^{\el}}:\ml\longrightarrow H^{\el}$$ 
\bigbreak
We begin by introducing an $l$-twisted version of the hypertoric variety.


\bigbreak
\subsection{The $l$-twisted multiplicative Hypertoric variety}
Where $q$ is a primitive $l^{th}$ root of unity, we have seen that $\qdiff^{\circ}$ is a sheaf of Azumaya algebras on $T^{*}\bA^{n,(l),\circ}$. From this, we can construct Azumaya algebras on $\fM^{\paren{l}}_{\chi,\eta}$, an $l$-twisted version of the classical multiplicative hypertoric variety. We first recall its definition.
\bigbreak
Recall the quantum moment map $\Phi_{q,K}:\cO\paren{K}\rightarrow\qdiff^{\circ}$ in $\text{Rep}_{q}\paren{K}$, 
\eqn{\Phi_{q,K}:\cO\paren{K}\overset{\phi^{t,*}}{\longrightarrow}\cO\paren{T}\overset{\Phi_{q,T}}{\longrightarrow}\qdiff^{\circ}}
\eqn{u_{j}\mapsto\prod_{i=1}^{n}\paren{1+x_{i}\partial_{i}}^{a_{ij}}}
corresponding to an embedding $\phi:K\hookrightarrow T$ given by $\phi\paren{k}_{i}=\prod^{d}_{j=1}k^{a_{ij}}_{j}$. From Lemma~\ref{powerlem} we can define the associated $l$-twisted moment map 
\eqn{\Phi^{\paren{l}}_{q,K}:\cO\paren{K^{(l)}}\longrightarrow Z_{l,\Delta^{l}}} by
\eqn{u^{l}_{j}\mapsto\prod^{n}_{i=1}\paren{1+x^{l}_{i}\partial^{l}_{i}}^{a_{ij}}} 
This corresponds to a morphism of varieties $\paren{T^{*}\bA^{n}}^{(l),\circ}\rightarrow K^{\paren{l}}$ and we have the following commutative diagram
\begin{center}
\begin{tikzpicture}[node distance=4.5cm, on grid]
\node (A) at (0,0) {$\cO\paren{K^{(l)}}$};
\node (B) [right=of A] {$Z_{l,\Delta^{l}}$};
\node (C) [below=of A,yshift=2.5cm] {$\cO\paren{K}$};
\node (D) [right=of C] {$\sO_{\paren{T^{*}\bA^{n}}^{(l),\circ},q}\paren{\paren{T^{*}\bA^{n}}^{(l),\circ}}$};
\draw[->] (A) -- node [above] {$\Phi^{\paren{l}}_{q,K}$} (B);
\draw[right hook->] (A) -- node [left] {} (C);
\draw[right hook->] (B) -- node [right] {} (D);
\draw[->] (C) -- node [below] {$\text{gr}(\Phi_{q,K})$} (D);
\end{tikzpicture}
\end{center}
We can now give the following
\defn{Let $K\subset T$ be a connected subtorus of $T$, $\chi\in X^{*}\paren{K}$ and $\eta\in K$. The $l$-twisted multiplicative hypertoric variety $\fM^{\paren{l}}_{\chi,\eta}$ is given by the GIT quotient 
\eqn{\fM^{\paren{l}}_{\chi,\eta}:=\paren{\Phi^{\paren{l}}_{K}}^{-1}\paren{\eta^{l}}\git_{\chi^{\paren{l}}} K^{\paren{l}}}}
Here, the character $\chi^{\paren{l}}$ of $K^{\paren{l}}$ is defined by $\chi^{\paren{l}}\paren{Fr_{l}(k)}=\chi\paren{k^{l}}$. As with the usual hypertoric variety, we can describe this $l$-twisted variety as
\eqn{\fM^{\paren{l}}_{\chi,\eta}=\text{Proj}\paren{\bigoplus_{m\geq 0}\paren{\sO\paren{X^{\paren{l}}}\slash\paren{\Phi^{\paren{l}}_{q,K}-\eta^{l}}}^{\paren{\chi^{\paren{l}}}^{m}}}}
The commutative diagram 
\begin{center}
\begin{tikzpicture}[node distance=3cm,on grid]
\node (A) at (0,0) {$\cO\paren{K^{\paren{l}}}$};
\node (B) [right=of A] {$\cO\paren{K}$};
\node (C) [below=of A,yshift=1.5cm] {$\sO\paren{\paren{T^{*}\bA^{n}}^{(l),\circ}}$};
\node (D) [right=of C] {$\sO\paren{\paren{T^{*}\bA^{n}}^{\circ}}$};
\draw[->] (A) -- (B);
\draw[->] (A) -- (C);
\draw[->] (B) -- (D);
\draw[->] (C) -- (D);
\end{tikzpicture}
\end{center}
induces a finite map of $\bK$-schemes $\text{Fr}_{l}:\fM_{\chi,\eta}\longrightarrow\fM^{\paren{l}}_{\chi,\eta}$. 
\subsection{The Azumaya algebra $\adiff$ on $\fM^{\paren{l}}_{\chi,\eta}$}
We have shown that $\qdiff^{\circ}$ is an Azumaya algebra over its centre $Z_{l,\Delta^{l}}$ and thus defines a locally free coherent sheaf over $Spec\paren{Z_{l,\Delta^{l}}}\cong  T^{*}\bA^{n,(l),\circ}$. By taking a quantum Hamiltonian reduction of this algebra, we can construct an Azumaya algebra on the $l$-twisted multiplicative hypertoric variety $\ml$. Let $\paren{K,\chi,\eta}$ be a triple as before and let $X:=\Phi^{-1}_{K}\paren{\eta}$, an affine variety, with algebra of functions $\sO\paren{X}=k\sqbrac{x_{i},\partial_{i}}\slash k\sqbrac{x_{i},\partial_{i}}(\Phi^{\#}_{K}-\eta)$. Let $X^{ss}$ be the semi-stable locus with respect to the character $\chi$. We thus have that $X^{(l)}=\Phi^{-1}_{K^{(l)}}\paren{\eta^{l}}$ and that $(X^{ss})^{(l)}=(X^{(l)})^{\chi^{(l)}-ss}$. Denote by $\varpi:X^{ss}\rightarrow\mh$ the geometric quotient morphism whose fibers are closed $K$-orbits. From this, we obtain the geometric quotient morphism $\varpi^{(l)}:(X^{ss})^{(l)}\rightarrow\ml$ with fibers as closed $K^{(l)}$-orbits. We make the following definition
\defn{The sheaf of $q$-difference operators on $\ml$ is given by \eqnn{\label{adiff}\adiff:=\varpi^{\el}_{*}\paren{\qdiff^{\circ}\paren{\bA^{n}}\slash\qdiff^{\circ}\paren{\bA^{n}}(\Phi_{q,K}-\eta)\vert_{\paren{X^{ss}}^{\el}}}^{K}}} Note that $\qdiff^{\circ}\paren{\bA^{n}}\slash\qdiff^{\circ}\paren{\bA^{n}}(\Phi_{q,K}-\eta)$ is an $\sO\paren{X^{\el}}$-module and hence a coherent sheaf on $X^{\el}$. As a sheaf on $\ml$ therefore, $\adiff$ has sections $$\Gamma\paren{U,\adiff}=\Gamma\paren{[\varpi^{\el}]^{-1}\paren{U},\qdiff^{\circ}\paren{\bA^{n}}\slash\qdiff^{\circ}\paren{\bA^{n}}(\Phi_{q,K}-\eta)}^{K}$$ for $U\subset\ml$ open. Moreover, we have the following
\prop{\label{sheafalg}$\adiff$ is a sheaf of algebras on $\ml$}
\proof{From the definition of $\ml$, one can see that it has an open affine cover by sets of the form $$U_{s}:=Spec\paren{\sO\paren{X^{\el}}\sqbrac{s^{-1}}^{K^{\paren{l}}}}$$ for $s\in\sO\paren{X^{\el}}^{\paren{\chi^{\paren{l}}}^{n}}$ for $n>0$. Choosing such an $s$ and a lift of it $\tilde{s}\in Z_{l,\Delta^{l}}$, we can compute the sections on $U_{s}$ as follows
\eqn{\Gamma\paren{U_{s},\adiff}&=\Gamma\paren{Spec\paren{\sO\paren{X^{\el}}\sqbrac{s^{-1}}},\qdiff^{\circ}}^{K}=\paren{\qdiff^{\circ}\otimes_{\sO\paren{X^{\el}}}\sO\paren{X^{\el}}\sqbrac{s^{-1}}}^{K}\\
&\cong\paren{\qdiff^{\circ}[\tilde{s}^{-1}]\slash\qdiff^{\circ}[\tilde{s}^{-1}]\paren{\Phi_{s}-\eta}}^{K}} Here, we are again considering $\qdiff^{\circ}$ as a locally free coherent sheaf on $X^{\paren{l}}$ and $\Phi_{s}$ is the map $$\cO\paren{K}\overset{\Phi_{q,K}}{\longrightarrow}\qdiff^{\circ}\rightarrow\qdiff^{\circ}[s^{-1}]$$ These sections thus have an algebra structure, from the general fact that a Hamiltonian reduction of an associative algebra inherits its algebra structure. As in Proposition~\ref{mh-H-mom} where the affine charts have non-empty intersection, this is given by $U_{r}\cap U_{s}=U_{rs}$ and the algebra structures agree when restricted to these intersections. Hence, these structures glue to give a sheaf of algebras on $\ml$.}
\bigbreak
The sheaf $\adiff$ possesses the following important property
\prop{\label{locfree}There is an isomorphism $\adiff\cong\paren{Fr_{l}}_{*}\sO_{\mathfrak{M}_{\chi,\eta}}$ as coherent $\sO_{\ml}$-modules.}
\proof{Consider the following commutative diagram
\begin{center}
\begin{tikzpicture}[node distance=3cm,on grid]
\node (A) at (0,0) {$X^{ss}$};
\node (B) [right=of A] {$\mh$};
\node (C) [below=of A,yshift=1cm] {$(X^{ss})^{\el}$};
\node (D) [right=of C] {$\ml$};
\draw[->] (A) -- node [above] {$\varpi$} (B);
\draw[->] (A) -- node [left]  {$Fr_{l}$} (C);
\draw[->] (B) -- node [right] {$Fr_{l}$} (D);
\draw[->] (C) -- node [below] {$\varpi^{\el}$} (D);
\end{tikzpicture}
\end{center}
where the top and bottom arrows are principal $K$- and $K^{\el}$-bundles respectively. Recall that $\sO\paren{X^{\el}}=Z_{l}\slash Z_{l}\cdot\paren{\Phi^{\#}_{K^{\paren{l}}}-\eta^{l}}$, and note that $\qdiff^{\circ}\slash\qdiff^{\circ}\cdot\paren{\Phi_{q,K}-\eta}$ is isomorphic to $(Fr_{l})_{*}\sO\paren{X}$ as both $\sO\paren{X^{\el}}$-modules and as $K$-modules. Now, given that $\sO_{\mh}=\varpi_{*}\paren{\sO\paren{X}\vert_{X^{ss}}}^{K}$, the commutativity of the diagram gives isomorphisms of $\sO_{\ml}$-modules:
\begin{equation*}
\begin{aligned}
\adiff &=\varpi^{\el}_{*}\paren{\qdiff^{\circ}\paren{\bA^{n}}\slash\qdiff^{\circ}\paren{\bA^{n}}(\Phi_{q,K}-\eta)\vert_{\paren{X^{ss}}^{\el}}}^{K} \\
       &\cong\varpi^{\el}_{*}\paren{\paren{Fr_{l}}_{*}\sO\paren{X}\vert_{\paren{X^{ss}}^{\el}}}^{K} \\
       &\cong\paren{Fr_{l}}_{*}\varpi_{*}\paren{\sO\paren{X}\vert_{X^{ss}}}^{K} \\
       &=\paren{Fr_{l}}_{*}\sO_{\mathfrak{M}_{\chi,\eta}} \\
\end{aligned}
\end{equation*}
}


\remark{The previous proposition says that $\adiff$ is a deformation of the $l$-twisted multiplicative hypertoric variety. In essentially following from the definitions, our proof is another illustration of how working with abelian groups significantly simplifies the situation. Proving an analogous result for non-abelian reductive groups would be considerably more involved as there is no counterpart of the initial isomorphism of $X^{*}(K)$-graded modules above. The positive characteristic analogue of this result for non-abelian groups, discussed and proved in 4.2 and 4.3 of \cite{BFG}, requires the formation of a \textit{Rees algebra} version $\mathbbm{\cR}\mathscr{A}_{\chi}$ of their Hamiltonian reduction sheaf $\mathscr{A}_{\chi}$. This is a deformation defined over $\bA^{1}\times\mathfrak{M}^{(1)}$, flat in the $\bA^{1}$-factor, and such that $\mathbbm{\cR}\mathscr{A}_{\chi}\vert_{\brac{0}\times\mathfrak{M}^{(1)}}\cong Fr_{*}\sO_{\mathfrak{M}}$ and $\mathbbm{\cR}\mathscr{A}_{\chi}\vert_{\brac{1}\times\mathfrak{M}^{(1)}}\cong\mathscr{A}_{\chi}$. They thus obtain that $Fr_{*}\sO_{\mathfrak{M}}$ and $\mathscr{A}_{\chi}$ represent the same class in the Grothendieck group. It is expected that adopting their approach should yield the appropriate analogue of Proposition~\ref{locfree} for non-abelian groups in our setting.}
\bigbreak
The rest of this section is devoted to proving 
\theo{\label{adiffaz}Suppose $\ml$ is smooth. Then $\adiff$ is a sheaf of Azumaya algebras on $\ml$.} 
Consider the short exact sequence of groups $$1\longrightarrow K_{l}\longrightarrow K\overset{Fr_{l}}{\longrightarrow} K^{(l)}\longrightarrow 1$$ where $K_{l}\cong\paren{\bZ\slash l\bZ}^{d}$. We will demonstrate that $\adiff$ is Azumaya by using the basic phenomenon, highlighted in $\cite{BFG}$, that a quantum Hamiltonian reduction with respect to $K$, as appears in the definition of $\adiff$, can also be achieved in two stages, first as a reduction with respect to $K_{l}$ and then via a reduction with respect to $K^{(l)}\cong K\slash K_{l}$. One can then reduce the question of $\adiff$ being Azumaya to that of the algebra arising from the $K_{l}$-reduction being Azumaya. We will now describe this two-stage reduction. 
\bigbreak
Recall that for $\eta\in K$ we can form the ideal $\ieta\trianglelefteq\cO\paren{K}$ generated by the elements $\brac{u_{j}-\eta_{j}}_{1\leq j\leq d}$. Consider now the ideal $\ietal:=\ieta\cap\cO\paren{K^{\paren{l}}}$. Setting $\rm{o}_{\eta}\paren{K}:=\cO\paren{K}\slash\cO\paren{K}\cdot\ietal$, we have a short exact sequence
\eqnn{\label{aexseq}0\longrightarrow\cO\paren{K}\cdot\ietal\longrightarrow\cO\paren{K}\overset{r}{\longrightarrow}\rm{o}_{\eta}\paren{K}\longrightarrow 0} Denote by $\sD^{\circ}_{q,\eta}:=\qdiff^{\circ}\slash\qdiff^{\circ}\cdot\Phi_{q,K}(\ietal)$ which is just the sheaf $\qdiff^{\circ}$ restricted to the closed affine subvariety $X^{\el}$. Given the short exact sequence~\ref{aexseq} the quantum moment map $\Phi_{q,K}:\cO\paren{K}\rightarrow\qdiff^{\circ}$ descends to a $K_{l}$-equivariant algebra homomorphim $\Phi_{K_{l}}:\rm{o}_{\eta}\paren{K}\rightarrow\sD^{\circ}_{q,\eta}$. Let $\mathfrak{i}_{\eta}:=r\paren{\ieta}\subset\rm{o}_{\eta}\paren{K}$. We can now form the Hamiltonian reduction with respect to the $K_{l}$-action, to obtain
$$\mathscr{E}_{\eta}:=\paren{\sD^{\circ}_{q,\eta}\slash\sD^{\circ}_{q,\eta}\cdot\Phi_{K_{l}}\paren{\mathfrak{i}_{\eta}}}^{K_{l}}$$ This is a sheaf of algebras on $X^{\el}$ which, given that we have restricted to $K_{l}$-invariants, has a $K^{\paren{l}}$-equivariant structure. We have the following version of Lemma 3.6.4. from \cite{BFG}: 
\prop{There is a canonical algebra isomorphism $$\Gamma\paren{T^{*}\bA^{n,\paren{l}},\qdiff^{\circ}\slash\qdiff^{\circ}\cdot\Phi_{q,K}\paren{\ieta}}^{K}\cong\Gamma\paren{X,\mathscr{E}_{\eta}}^{K^{\paren{l}}}$$}
\proof{Observe that there is an isomorphism between $$\qdiff^{\circ}\cdot\Phi_{q,K}\paren{\ieta}\slash\qdiff^{\circ}\cdot\Phi_{q,K}\paren{\ietal}\cong\sD^{\circ}_{q,\eta}\cdot\Phi_{K_{l}}\paren{\mathfrak{i}_{\eta}}$$ and so we have that $$\qdiff^{\circ}\slash\qdiff^{\circ}\cdot\Phi_{q,K}\paren{\ieta}\cong\sD^{\circ}_{q,\eta}\slash\sD^{\circ}_{q,\eta}\cdot\Phi_{K_{l}}\paren{\mathfrak{i}_{\eta}}$$ 
Taking $K_{l}$-invariants thus gives $$\paren{\qdiff^{\circ}\slash\qdiff^{\circ}\cdot\Phi_{q,K}\paren{\ieta}}^{K_{l}}\cong\mathscr{E}_{\eta}$$ Hence, we have that 
\begin{equation*}
\begin{aligned}
\Gamma\paren{T^{*}\bA^{n,\paren{l}},\qdiff^{\circ}\slash\qdiff^{\circ}\cdot\Phi_{q,K}\paren{\ieta}}^{K} & \cong\Gamma\paren{T^{*}\bA^{n,\paren{l}},\paren{\qdiff^{\circ}\slash\qdiff^{\circ}\cdot\Phi_{q,K}\paren{\ieta}}^{K_{l}}}^{K^{\paren{l}}} \\
& \cong\Gamma\paren{T^{*}\bA^{n,\paren{l}},\mathscr{E}_{\eta}}^{K^{\paren{l}}} \\
& \cong\Gamma\paren{X,\mathscr{E}_{\eta}}^{K^{\paren{l}}} \\
\end{aligned}
\end{equation*}
}
From the above computation, we can see that $\adiff\cong\varpi^{\el}_{*}\paren{\mathscr{E}_{\eta}\vert_{(X^{ss})^{\el}}}^{K^{\paren{l}}}$. Proposition~\ref{locfree} implies that $\adiff$ is locally free as an $\sO_{\ml}$-module.\\
Since $\varpi^{\el}:(X^{ss})^{\el}\rightarrow\ml$ is a geometric quotient morphism, we have an isomorphism $[\varpi^{\el}]^{*}\adiff\cong\mathscr{E}_{\eta}\vert_{(X^{ss})^{\el}}$, and thus $\mathscr{E}_{\eta}\vert_{(X^{ss})^{\el}}$ is also locally free. Let $x\in (X^{ss})^{\el}$ and let $\bar{x}=\varpi^{\el}\paren{x}$. From the isomorphism above, there is an algebra isomorphism of fibers $\mathscr{E}_{x}\cong\mathscr{A}_{\bar{x}}$. Thus, in order to prove that $\adiff$ is an Azumaya algebra, it suffices to show that the fiber $\mathscr{E}_{x}$ is a matrix algebra for all $x\in (X^{ss})^{\el}$.
\bigbreak
We know that $\qdiff$ is generically Azumaya on $T^{*}\bA^{n,\paren{l}}$, so there is a finite dimensional vector space $V$ such that $\sD^{\circ}_{q,x}\cong\text{End}_{k}V$ for $x\in (X^{ss})^{\el}$. By definition, we have $$\mathscr{E}_{x}=\paren{\sD^{\circ}_{q,x}\slash\sD^{\circ}_{q,x}\cdot\mathfrak{i}_{\eta}}^{K_{l}}\cong\paren{\text{End}_{k}V\slash\text{End}_{k}V\cdot\mathfrak{i}_{\eta}}^{K_{l}}$$
To finish the proof, it suffices to show that $$\paren{\text{End}_{k}V\slash\text{End}_{k}V\cdot\mathfrak{i}_{\eta}}^{K_{l}}\cong\text{End}_{k}\paren{V_{\eta}}$$ where $V_{\eta}:=\brac{x\in V\vert u_{j}\paren{x}=\eta_{j}\paren{x}, 1\leq j\leq d}$ is the $\eta$-weight space of $V$.
\bigbreak
For ease of notation, let $\bD:=\sD^{\circ}_{q,x}\cong\text{End}_{k}V$ and let $\bD_{\eta}:=\paren{\bD\slash\bD\cdot\mathfrak{i}_{\eta}}^{K_{l}}$, its quantum Hamiltonian reduction with respect to $\eta\in K$ and the $K_{l}$-equivariant map $\Phi_{K_{l}}:\rm{o}_{\eta}\paren{K}\rightarrow\sD_{q,\eta}$. The action of $\bD$ on $V$ descends to a well-defined $\bD_{\eta}$-action on $V_{\eta}$. In addition, there is a right $\bD_{\eta}$-action on $\bD\slash\bD\cdot\mathfrak{i}_{\eta}$. We have the following
\lem{\label{res} The restriction map $\textbf{res}:\text{End}_{k}V\rightarrow\text{Hom}_{k}\paren{V_{\eta},V}$ induces an isomorphism of $\bD-\bD_{\eta}$-bimodules $$\bD\slash\bD\cdot\mathfrak{i}_{\eta}\overset{\sim}{\longrightarrow}\text{Hom}_{k}\paren{V_{\eta},V}$$}
\proof{There is a correspondence between vector subspaces $U\subset V$ and left ideals $J\triangleleft\bD$ as follows: to a subspace $U\subset V$ let $J_{U}$ be the left ideal given by $\brac{f\in\text{End}_{k}V\vert f\vert_{U}=0}$; to any left ideal $J\triangleleft\bD$, one has that $J=J_{U}$ where $U$ is the subspace $U=\bigcap_{f\in J}\ker{f}$. If we apply this correspondence to $J=\bD\cdot\mathfrak{i}_{\eta}$, we have that $J=J_{U}$ where $U=\bigcap_{x\in\mathfrak{i}_{\eta}}\ker{\Phi_{K_{l}}\paren{x}}$ which, by definition, is the weight space $V_{\eta}$. Hence, we have that $\bD\cdot\mathfrak{i}_{\eta}=\text{Hom}_{k}\paren{V\slash V_{\eta},V}$ and so $\bD\slash\bD\cdot\mathfrak{i}_{\eta}\cong\text{Hom}_{k}\paren{V_{\eta},V}$.} Given this isomorphism, we can now show
\lem{\label{lriso}There are isomorphisms of algebras $$\bD_{\eta}\overset{\overset{\textbf{res}}{\sim}}{\longrightarrow}\text{End}_{\bD}\paren{\text{Hom}_{k}\paren{V_{\eta},V}}^{opp}\overset{\sim}{\longleftarrow}\text{End}_{k}V_{\eta}$$}
\proof{Let $u\in\bD$ be such that $u\paren{\text{mod}\paren{\bD\cdot\mathfrak{i}_{\eta}}}\in\bD_{\eta}$. The map $$u\mapsto\paren{f_{u}:\bD\slash\bD\cdot\mathfrak{i}_{\eta}\rightarrow\bD\slash\bD\cdot\mathfrak{i}_{\eta}}$$ where $f_{u}$ denotes right multiplication by $u$, induces an algebra isomorphism $\bD_{\eta}\rightarrow\paren{\text{End}_{\bD}\paren{\bD\slash\bD\cdot\mathfrak{i}_{\eta}}}^{opp}$. Combining this with Lemma~\ref{res} we obtain the first isomorphism above. The second isomorphism is given by the map $$\text{End}_{k}V_{\eta}\rightarrow\text{End}_{\bD}\paren{\text{Hom}_{k}\paren{V_{\eta},V}}^{opp}$$ defined by $f\mapsto\paren{g\mapsto g\circ f}$.} \bigbreak
Combining these two lemmas we have that $\mathscr{E}_{x}$ is a matrix algebra and hence, that $\adiff$ is an Azumaya algebra on $\ml$.
\subsection{Azumaya splitting on fibers}
Having constructed an Azumaya algebra $\adiff$, on the $l$-twisted hypertoric variety $\ml$, we will now study its splitting behaviour. Recall that an Azumaya algebra $\sA$ on a scheme $X$ is called \textit{split} if there is a vector bundle $\sV$ on $X$ and an $\sO_{X}$-module isomorphism $\sA\cong\mathscr{E}nd_{\sO_{X}}\paren{\sV}$. We will demonstrate that $\adiff$ splits when restricted to the fibers of the resolution map $$\Psi:\ml\longrightarrow\Phi_{K^{\el}}^{-1}\paren{\eta^{l}}\slash K^{\el}$$ and the moment map $$\Phi_{H^{\el}}:\ml\longrightarrow H^{\el}$$ To do this, we adapt the approach used in \cite{Sta1} in which, working over a field of prime characteristic, an \'etale cover of the Frobenius twisted hypertoric variety is constructed, splitting the analogue of our sheaf $\adiff$. That this gives a splitting of the Azumaya algebra on fibers of the resolution map is implicit in \cite{Sta1} - a fact which, together with its adaptation to our present setting, was outlined to me in correspondence with Iordan Ganev, and appears in his preprint \cite{Ga}. We begin with 
\defn{Let $\sA$ be an Azumaya algebra on a scheme $Y$. An \'etale splitting of $\sA$ is a morphism of schemes $f:X\rightarrow Y$ such that $f^{*}\sA$ is split on $X$.} Denote by $R:=Z_{l,\Delta^{l}}\sqbrac{\brac{\alpha_{i}}}$ the commutative, though non-central, sub-algebra of $\qdiff^{\circ}$, obtained by adjoining the Euler operators to $Z_{l,\Delta^{l}}$. We have the following result which points towards what we expect to be the case for the Hamiltonian reduction sheaf. It can be proved in a similar manner to the analogous theorem 3.2.5 in \cite{Sta1}.
\prop{\label{d-etale}The covering map $f:T^{*}\bA^{n,\el,\circ}\times_{T^{\el}}T\rightarrow T^{*}\bA^{n,\el,\circ}$ induced by the inclusion $Z_{l,\Delta^{l}}\hookrightarrow R$ is an \'etale splitting of $\qdiff^{\circ}$.} 
\rem{The key observation in \cite{Sta1} is that the Euler operators $$\brac{x_{i}\partial_{i}\vert 1\leq i\leq n}\subset\sD\paren{\bA^{n}}$$ satisfy the Artin-Schreier equation - that is, $\paren{x_{i}\partial_{i}}^{p}-x_{i}\partial_{i}-x_{i}^{p}\partial_{i}^{p}=0$ - with coefficients in $\sO_{T^{*}\bA^{n,\el}}$. This allows the $T$-equivariant \'etale splitting cover - analogous to that in Proposition~\ref{d-etale} - constructed therein to descend to such a cover of the Frobenius twisted hypertoric variety. This property, and hence such an \'etale cover, does not have an analogue for non-abelian reductive groups.} 
The role of the Artin-Schreier map in our context is played by the $l^{th}$ power map $Fr_{l}$ and we have, from Lemma~\ref{delta-l}, that the $\alpha_{i}$'s satisfy $\alpha_{i}^{l}=\paren{1+x_{i}\partial_{i}}^{l}=1+x_{i}^{l}\partial_{i}^{l}$. As such, we can similarly construct an \'etale cover of $\ml$, splitting $\adiff$ and inducing, via Lemma~\ref{Cart-Az}, splittings of its restriction to the maps mentioned at the beginning of this section. 
\bigbreak
\defn{\label{mathscrY}Define $\mathscr{Y}_{\chi,\eta}:=\ml\times_{H^{\el}}H$.}
The factor product $\ycov$ fits into the Cartesian square
\begin{center}
\begin{tikzpicture}[node distance=2cm,on grid]
\node (A) at (0,0) {$\ycov$};
\node (B) [right=of A] {$\ml$};
\node (C) [below=of A,yshift=0.5cm] {$H$};
\node (D) [right=of C] {$H^{\el}$};
\draw[->] (A) -- node [above] {$f_{\mathscr{Y}}$} (B);
\draw[->] (A) -- node [left]  {$\Phi_{R,H}$} (C);
\draw[->] (B) -- node [right] {$\Phi_{H^{\el}}$} (D);
\draw[->] (C) -- node [below] {$Fr_{l}$} (D);
\end{tikzpicture}
\end{center}
and we see immediately that the covering map $f_{\mathscr{Y}}$ is \'etale of degree $l^{n-d}$.
Denote by $\Phi_{R,K}:Spec\paren{R}\rightarrow K$ the composition $Spec\paren{R}\overset{\Phi_{R,T}}{\rightarrow}T\overset{\phi^{t}}{\rightarrow}K$ and let $Y:=\Phi^{-1}_{R,K}\paren{\eta}$. We observe the following: 
\prop{\label{GITiso} The GIT quotient $Y\git_{\chi^{\el}}K^{\el}$ is isomorphic to $\ycov$.} 
\proof{Restricting the map induced by the injection $Z_{l,\Delta^{l}}\hookrightarrow R$ gives a map $f:Y\rightarrow X^{\el}$. We have the following explicit descriptions of the global sections of these spaces:
\eqn{\sO\paren{X^{\el}}\cong Z_{l,\Delta^{l}}\slash\paren{\Phi_{K^{\el}}^{\#}-\eta^{l}}} and
\eqnn{\label{O(Y)O(X)}\sO\paren{Y}\cong\sO\paren{X^{\el}}\sqbrac{\brac{T_{i}^{\pm}}}\slash\paren{T_{i}^{l}-\alpha_{i}^{l},\prod_{i=1}^{n}T_{i}^{a_{ij}}-\eta_{j}}}
From these formulae, we can readily see that $Y\cong X^{\el}\times_{H^{\el}}H$ with respect to the group-valued moment maps $Y\rightarrow H$ and $X^{\el}\rightarrow H^{\el}$. To complete the proof, let $s\in\sO\paren{X^{\el}}^{\paren{\chi^{\el}}^{n}}$ for some $n>0$.  Note that, as the $K^{\el}$-action on $Spec\paren{R}$ is induced by that on $Spec\paren{Z_{l}}$, there is a surjection $Y\git_{\chi^{\el}}K^{\el}\longrightarrow\ml$ given locally by $$Spec\paren{\sO\paren{Y}\sqbrac{s^{-1}}^{K^{\el}}}\longrightarrow Spec\paren{\sO\paren{X^{\el}}\sqbrac{s^{-1}}^{K^{\el}}}$$ Again, by the explicit formulae for the global sections algebras given above, we are done.}
\bigbreak
Let $\varpi_{\sY}:Y^{ss}\rightarrow\ycov$ denote the quotient map and define a sheaf on $\ycov$ by $\sD_{q,\sY}:=\varpi_{\sY*}\paren{\sD^{\circ}_{q}\vert_{Y^{ss}}}^{K}$.
\lem{The sheaf $\sD_{q,\sY}$ is locally free on $\ycov$.}
\proof{Since $\sD_{q,\sY}$ is a coherent sheaf on a scheme of finite type over $\bK$, it suffices to show that the dimension of the fiber of $\sD_{q,\sY}$ is equal at all closed points of $\ycov$. Let $p\in\ycov$ be a closed point lying under $p'\in Y^{ss}$. Since $\ml$ is smooth, $K^{\el}$ acts freely on the orbit of $p'$ in $Y^{ss}$. One can identify the fiber $(\sD_{q,\sY})_{p}$ with $(\qdiff^{\circ}\otimes_{R}k(p'))^{K_{l}}$, the $K_{l}$-invariant sections of the fiber of $\qdiff^{\circ}$ over $p'$. One can then show that the latter space is of dimension $l^{n-d}$ over $k$.}
The main result of this section is the following:
\theo{The map $f_{\mathscr{Y}}:\ycov\rightarrow\ml$ is an \'etale splitting of the Azumaya algebra $\adiff$ via an isomorphism $(f_{\sY})^{*}\adiff\overset{\sim}{\longrightarrow}\mathscr{E}nd_{\ycov}\paren{\sD_{q,\sY}}$.}
\proof{We use the GIT description from Proposition~\ref{GITiso}. \\
We construct a map $(f_{\sY})^{*}\adiff\longrightarrow\mathscr{E}nd_{\ycov}\paren{\sD_{q,\sY}}$ as follows.  Let $s\in\sO\paren{X^{\el}}^{\paren{\chi^{\el}}^{n}}$ and set
\eqn{U_{s}=Spec\paren{\sO\paren{X^{\el}}\sqbrac{s^{-1}}^{K^{\el}}}} and
\eqn{V_{s}=Spec\paren{\sO\paren{Y}\sqbrac{s^{-1}}^{K^{\el}}}}
From the previous proposition, we have that $\ycov$ is covered by open affine sets of the form $V_{s}$ and that $f_{\sY}^{-1}\paren{U_{s}}=V_{s}$. From the isomorphism of $\sO(X^{\el})$-algebras, $\qdiffc\slash(\Phi_{q,K}-\eta)\cong\qdiffc\otimes_{R}\sO(Y)$ we can form an isomorphism of $\sO(U_{s})$-algebras 
$$\paren{\qdiffc\slash(\Phi_{q,K}-\eta)\otimes_{\sO(X^{\el})}\sO(X^{\el})\sqbrac{s^{-1}}}^{K}\cong\paren{\qdiffc\otimes_{R}\sO(Y)\sqbrac{s^{-1}}}^{K}$$ and hence a map of $\sO(U_{s})$-algebras $$\paren{\qdiffc\slash(\Phi_{q,K}-\eta)\otimes_{\sO(X^{\el})}\sO(X^{\el})\sqbrac{s^{-1}}}^{K}\longrightarrow\text{End}_{\sO(V_{s})}\paren{\qdiffc\otimes_{R}\sO(Y)\sqbrac{s^{-1}}^{K}}$$
Finally, tensoring by $\sO(V_{s})$ gives a map of $\sO(V_{s})$-algebras 
$$\paren{\qdiffc\slash(\Phi_{q,K}-\eta)\otimes_{\sO(X^{\el})}\sO(X^{\el})\sqbrac{s^{-1}}}^{K}\otimes_{\sO(U_{s})}\sO(V_{s})\longrightarrow\text{End}_{\sO(V_{s})}\paren{\qdiffc\otimes_{R}\sO(Y)\sqbrac{s^{-1}}^{K}}$$
Since $f_{\sY}$ is surjective, $f_{\sY}(V_{s})=U_{s}$. Recalling the definition of $\adiff$, the source of the above map gives sections of the sheaf $(f_{\sY})^{*}\adiff$ over $V_{s}$. The target of the map describes sections of the sheaf $\mathscr{E}nd_{\ycov}\paren{\sD_{q,\sY}}$ over $V_{s}$. This map is compatible with gluing along members of the affine cover of $\ycov$, and hence we obtain a morphism of sheaves on $\ycov$. \bigbreak
It remains to show that this map is an isomorphism and, since we are mapping between locally free coherent sheaves, we can do this locally. Again, let $p'\in Y^{ss}$ be a closed point lying over $p\in\ycov$. The map induced on fibers by the map above is $(\sD^{\circ}_{q,p})^{K_{l}}\rightarrow End_{k}((\sD^{\circ}_{q,p'})^{K_{l}})$, which one can show is an isomorphism.}
\prop{The \'etale cover $f_{\mathscr{Y}}:\mathscr{Y}_{\chi,\eta}\rightarrow\ml$ fits into the Cartesian square:
\begin{center}
\begin{tikzpicture}[node distance=3.5cm,on grid]
\node (A) at (0,0) {$\ycov$};
\node (B) [right=of A] {$\ml$};
\node (C) [below=of A,yshift=2cm] {$\mathscr{Y}_{\chi,\eta, aff}$};
\node (D) [right=of C] {$\Phi_{K^{\el}}^{-1}\paren{\eta^{l}}\slash K^{\el}$};
\draw[->] (A) -- node [above] {$f_{\mathscr{Y}}$} (B);
\draw[->] (A) -- node [left]  {} (C);
\draw[->] (B) -- node [right] {$\Psi$} (D);
\draw[->] (C) -- node [below] {$f_{aff}$} (D);
\end{tikzpicture}
\end{center}}
\proof{Given the GIT description of $\ycov$, the proof is similar to that of Proposition~\ref{GITiso}. Once again, from equation~\ref{O(Y)O(X)}, we can see that $Y\cong X^{\el}\times_{X^{\el}\slash K^{\el}}Y\slash K^{\el}$ and, locally, we have the Cartesian square:
\begin{center}
\begin{tikzpicture}[node distance=5.5cm,on grid]
\node (A) at (0,0) {$Spec\paren{\sO\paren{Y}\sqbrac{s^{-1}}^{K^{\el}}}$};
\node (B) [right=of A] {$Spec\paren{\sO\paren{X^{\el}}\sqbrac{s^{-1}}^{K^{\el}}}$};
\node (C) [below=of A,yshift=3.5cm] {$Spec\paren{\sO\paren{Y}^{K^{\el}}}$};
\node (D) [right=of C] {$Spec\paren{\sO\paren{X^{\el}}^{K^{\el}}}$};
\draw[->] (A) -- node [above] {$f_{\mathscr{Y},s}$} (B);
\draw[->] (A) -- node [left]  {} (C);
\draw[->] (B) -- node [right] {} (D);
\draw[->] (C) -- node [below] {$f_{\mathscr{Y},aff}$} (D);
\end{tikzpicture}
\end{center}
and so we are done.}
\rem{Similar to \cite{Sta1}, we can fit the \'etale cover into the above Cartesian square because the functions on $Y$ are obtained by adjoining $K^{\el}$-invariant elements satisfying the equation $T_{i}^{l}-(1+x_{i}^{l}\partial_{i}^{l})=0$ to functions on $X^{\el}$. As such, the Cartesian property is preserved through the GIT reduction.}
Finally, via Lemma~\ref{Cart-Az}, we have the following:
\cor{\label{split-fiber} The Azumaya algebra $\adiff$ splits on the fibers of the moment map $\Phi_{H^{\el}}:\ml\longrightarrow H^{\el}$ and the resolution $\Psi:\ml\longrightarrow\Phi_{K^{\el}}^{-1}\paren{\eta^{l}}\slash K^{\el}$.}
Standard Morita Theory arguments then give the following:
\cor{\label{ab-eq-split}Let $\zeta\in\mlaff$ and $h\in H^{\el}$. There are equivalences of abelian categories
\eqn{\adiff\vert_{\Psi^{-1}(\zeta)}-\operatorname{mod}\simeq\sO_{\Psi^{-1}(\zeta)}-\operatorname{mod}\qquad\text{and }\qquad \adiff\vert_{\Phi_{H^{\el}}^{-1}(h)}-\operatorname{mod}\simeq\sO_{\Phi_{H^{\el}}^{-1}(h)}-\operatorname{mod}
} of coherent sheaves.}

\section{Localization for Quantum Hypertoric Varieties} \label{sec:loc}
In \cite{BeKu}, the authors describe a \textit{microlocalization} for central reductions of hypertoric enveloping algebras - that is, an equivalence of abelian categories between the category of finitely generated modules for the algebras introduced in section~\ref{subsection:HEA} and the category of coherent sheaves of $\mathscr{W}$-algebras on a hypertoric variety. In our setting we cannot obtain an equivalence at the level of abelian categories, however, as in \cite{BMR} and \cite{BFG}, treating the cases of Lie algebras and Cherednik algebras in positive characteristic respectively, and \cite{BK2} considering the case of quantum groups at roots of unity, we can obtain an equivalence of derived categories. In fact, the result we present here is a formal consequence of the techniques developed for proving such equivalences in \cite{BeKa} and \cite{BMR}, and so we present only an outline of the proof, directing the reader to the relevant sections of those papers for a full exposition. 
\bigbreak
Assume henceforth that $\ml$ is smooth and let $\cU_{q,\eta}:=\Gamma\paren{\ml,\adiff}$. Assume further that the affinization map $\Phi:\mh\rightarrow\mhaff$ is birational. Since $\ml$ is a finite-dimensional smooth variety $\Gamma:\adiff\operatorname{-mod}\rightarrow\cU_{q,\eta}\operatorname{-mod}$ has finite homological dimension, and so the derived functor $R\Gamma$ induces a map of the respective bounded derived categories. Denote by $\mathtt{Loc}:\cU_{q,\eta}\operatorname{-mod}\rightarrow\adiff\operatorname{-mod}$ the left adjoint functor to $\Gamma$ given by $$M\mapsto\adiff\otimes_{\underline{\cU_{q,\eta}}}\underline{M}$$ where $\underline{M}$ denotes the constant sheaf with stalk $M$, and let $$\mathscr{L}:D^{b}\paren{\cU_{q,\eta}\operatorname{-mod}}\rightarrow D\paren{\adiff\operatorname{-mod}}$$ be its left derived functor. 
\prop{The sheaf $\adiff$ has no higher global sections.}
\proof{Note that, since the affinization map $\Psi:\mh\rightarrow\mhaff$ is birational, the relative Grauert-Riemenschneider Theorem gives that $R^{i}\Psi_{*}\omega_{\mh}=0$ for $i>0$ and, since $\mh$ is a symplectic variety, $\omega_{\mh}\cong\sO_{\mh}$. We have the following sequence of isomorphisms of complexes of coherent sheaves of $\sO_{\ml}$-modules:
\begin{align*}
R\Gamma\paren{\ml,\adiff} 
&\overset{\text{Prop.}~\ref{locfree}}{\cong}R\Gamma\paren{\ml,(Fr_{l})_{*}\sO_{\mh}}		\cong R\Gamma\paren{\mh,\sO_{\mh}}\\ 
&\cong R\Gamma\paren{\mhaff,\Psi_{*}\sO_{\mh}	}														    \cong\Gamma\paren{\mhaff, \Psi_{*}\sO_{\mh}} \\
&\cong\Gamma\paren{\mh,\sO_{\mh}}																				\cong\Gamma\paren{\ml,(Fr_{l})_{*}\sO_{\ml}} \\
&\cong\Gamma\paren{\ml,\adiff} \\
\end{align*}}
\bigbreak
\remark{Here, we have assumed that we have chosen GIT parameters for which the affinization map $\Phi$ is birational and for which $\mh$ is smooth and, therefore, such that $\Phi$ is a symplectic resolution. In the additive setting, Bellamy and Kuwabara showed that the categorical quotient $\mathcal{M}_{1,0}$ is of dimension $2\paren{n-d}$ and is Cohen-Macaulay and that, for stability parameter $\chi\in X^{*}_{\mathbb{Q}}\paren{K}:=X^{*}\paren{K}\otimes_{\bZ}\mathbb{Q}$ lying in a maximal cone of the GIT fan, if the $\bZ$-valued matrix encoding the action of $K$ on $T^{*}\bA^{n}$ is unimodular, then the affinization map $f:\mathcal{M}_{\chi,0}\rightarrow\mathcal{M}_{1,0}$ is birational (see Proposition 4.6 and Lemma 4.9 of \cite{BeKu}). It is hoped that a similar argument to theirs can be applied in the multiplicative setting.}
 
\bigbreak
The context in which we have been working fits into an important paradigm in geometric representation theory by which, given an Azumaya algebra on a Calabi-Yau variety, one can obtain a derived equivalence between categories of coherent sheaves of modules over the Azumaya algebra and finitely generated modules over its global sections algebra. This formalism, involving a relative notion of Serre duality and Calabi-Yau triangulated category, is detailed in \cite{BeKa}, \cite{BMR}, and \cite{BMR2}. We now give a summary of these ideas before applying it to our situation. 
\bigbreak
Let $\sO$ be a finite type commutative algebra over $\bK$ and let $D$ be an $\sO$-linear triangulated category. We say that $D$ is an $\sO$-triangulated category if there is a functor $RHom_{D\slash\sO}:D^{op}\times D\rightarrow D^{b}\paren{mod^{fg}\paren{\sO}}$ and a functorial isomorphism $Hom_{D}\paren{X,Y}\cong H^{0}\paren{RHom_{D\slash\sO}\paren{X,Y}}$. An $\sO$-Serre functor on $D$ is an auto-equivalence $S:D\rightarrow D$ such that $$\mathbb{D}_{\sO}\paren{Hom_{D\slash\sO}\paren{X,Y}}\cong RHom_{D\slash\sO}\paren{Y,SX}$$ where $D_{\sO}$ is Grothendieck-Serre daulity for $\sO$-modules. An $\sO$-Calabi-Yau category is one for which the shift functor $X\rightarrow X\sqbrac{n}$ is an $\sO$-Serre functor for some $n\in\bZ$.
\bigbreak
The importance of Azumaya algebras is borne out by the following (Lemma 2.7 \cite{BeKa}):
\theo{Let $F:C\rightarrow D$ be a triangulated functor. Suppose that:
\begin{enumerate}
\item $F$ has a left adjoint $G$ such that $id\rightarrow F\circ G$ is an isomorphism.
\item $C$ is indecomposable in that it cannot be written as $C=C_{1}\oplus C_{2}$ for non-zero triangulated categories $C_{1}$ and $C_{2}$.
\item $C$ is an $\sO$-Calabi-Yau category for some finite-type, commutative algebra $\sO$
\end{enumerate}
then $F$ is an equivalence of categories.}
From \cite{BMR} Lemma 3.5.3 we have the following:
\prop{Let $\mathscr{A}$ be an Azumaya algebra over a connected, quasi-projective scheme $Y$ over a field $\bK$. The category $D^{b}\paren{\mathscr{A}-mod}$ is indecomposable.}
and from \cite{BMR2} Lemma 1.8.1. b.) we have:
\prop{Let $\sO$ be a finite-type commutative $\bK$-algebra. Let $Y$ be a $\bK$-variety with a projective morphism $\Psi:Y\rightarrow Spec\paren{\sO}$. Then $D^{b}\paren{\sO-mod}$ is an $\sO$-triangulated category by $RHom_{D\slash\sO}\paren{\mathscr{F},\mathscr{G}}:=R\Psi_{*}RHom(\mathscr{F},\mathscr{G})$. If $Y$ is smooth and quasi-projective, then
for an Azumaya algebra $\mathscr{A}$ on $Y$ , $D^{b}\paren{\mathscr{A}-mod}$ is an $\sO$-
triangulated category. The functor $\mathscr{F}\rightarrow \mathscr{F}\otimes\omega_{Y}\sqbrac{\operatorname{dim}Y}$ is an
$\sO$-Serre functor. If $Y$ is a Calabi-Yau variety then the $\sO$-triangulated category $D^{b}\paren{\mathscr{A}-mod}$ is Calabi-Yau.}
\remark{The key feature of Azumaya algebras which makes the proof of the above proposition relatively straightforward is that, after an \'etale base change, the Azumaya algebra is globally a matrix algebra and, therefore, has a trace form. This trace form descends to the Azumaya algebra over the original base. This allows one to define a non-degenerate pairing of $Hom$ spaces which in turn provides isomorphisms $$R\mathcal{H}om_{\mathscr{A}}\paren{\mathscr{G},\mathscr{F}}\cong R\mathcal{H}om_{\sO_{Y}}\paren{R\mathcal{H}om_{\mathscr{A}}\paren{\mathscr{F},\mathscr{G}},\sO_{Y}}$$ for any $\mathscr{F},\mathscr{G}\in D^{b}\paren{\mathscr{A}-mod}$ (see also \cite{BMR2} Lemma 1.8.1. c.)).}
\bigbreak
Taking $Y$ to be $\fM^{\paren{l}}_{\chi,\lambda}$ and $\Psi:\fM^{\paren{l}}_{\chi,\lambda}\rightarrow\Phi_{K^{\el}}^{-1}\paren{\eta^{l}}\slash K^{\el}$ to be the affinization morphism, and given that $\fM^{\paren{l}}_{\chi,\lambda}$ is a smooth symplectic variety and hence Calabi-Yau, we have:
\theo{\label{localization} If the functor $\mathtt{Loc}$ has finite homological dimension, the functors $\mathscr{L}$ and $R\Gamma$ give inverse equivalences of bounded derived categories

\begin{center}
\begin{tikzpicture}[node distance=5cm,on grid]
\path          node (A) {$D^{b}\paren{\adiff\operatorname{-mod}}$}
    (0:5cm)    node (G) {$D^{b}\paren{\cU_{q,\eta}\operatorname{-mod}}$};
 \path[-stealth]
   
   ([yshift=-2.5pt]A.east) edge node [above,yshift= 1.0ex]  {$R\Gamma$} ([yshift=-2.5pt]G.west)        
   ([yshift= 2.5pt]G.west) edge node [below,yshift=-1.0ex]  {$\mathscr{L}$} ([yshift= 2.5pt]A.east);
\end{tikzpicture}
\end{center}}
\rem{Throughout, we have restricted ourselves to parameters for which $\ml$ is smooth and hence such that the global sections functor has finite homological dimension. This is true generically, though it remains to provide an explicit characterization of this space of parameters.}
The algebra $\cU_{q,\eta}$ is a central reduction of a multiplicative $q$-deformed version of the hypertoric enveloping algebra at the character $\eta$. Classically, for certain parameters, there are isomorphisms between these Hamiltonian reduction algebras and familiar algebras of representation theoretic interest. As described in section 6 of \cite{BeKu}, in the additive setting, if we take $K$ to be an $(n-1)$-torus whose embedding into $T$ is given by the $(n-1)\times n$ matrix:
$$\begin{pmatrix}
1 & 0 & \ldots & 0 & -1 \\
 0 & 1 &   & \vdots & -1 \\
 \vdots & & \ddots & 0 & \vdots \\
 0 & \ldots & 0 & 1 & -1  
\end{pmatrix}$$
then, for $\chi\in X^{*}\paren{\mathfrak{k}}$ generic, away from certain GIT walls, $\cU_{q,\eta}$ is the spherical subalgebra of the rational Cherednik algebra associated to a cyclic group, which is generically Morita equivalent to the Cherednik algebra itself. The associated hypertoric variety is toric, and is a minimal resolution of the type A Kleinian singularity $\bK^{2}\slash\paren{\bZ\slash n\bZ}$. This can be described as the quiver variety associated to a cyclic quiver with $n$ vertices and dimension vector $d_{v}=1$ for all vertices $v$. In our setting, we obtain a $q$-deformed multiplicative version of this Cherednik algebra, whose representation theory we expect to be analogous to that of its positive characteristic cousin, as fully computed in \cite{La}.
\rem{In the case where $\ml$ is a minimal resolution of a type A Kleinian singularity, one could prove a version of the Azumaya splitting on fibers of the resolution from Corollary~\ref{split-fiber}, analogous to the proof given for the Hilbert-Chow map in Theorem 7.4.1. of \cite{BFG}, based on results in \cite{BeKa}. This provides an Azumaya splitting on formal neighbourhoods of the fibers of the resolution, $\widehat{\Psi^{-1}\paren{\zeta}}$. Combining the resulting equivalence of abelian categories \eqn{\adiff\vert_{\widehat{\Psi^{-1}(\zeta)}}-\operatorname{mod}\simeq\sO_{\widehat{\Psi^{-1}(\zeta)}}-\operatorname{mod}} from Corollary~\ref{ab-eq-split} with the derived localization, one obtains an equivalence of bounded derived categories 
\eqn{D^{b}\paren{\sO_{\widehat{\Psi^{-1}(\zeta)}}-\operatorname{mod}}\simeq D^{b}\paren{\widehat{\cU_{\eta, \zeta}}-\operatorname{mod}}}
of $\sO$-modules on these formal neighbourhoods, and finitely generated modules for $\widehat{\cU_{\eta, \zeta}}$, the completion of $\cU_{q,\eta}$ at the maximal ideal $\mathfrak{m}_{\zeta}$ of the $l$-centre, corresponding to the point $\zeta\in \sqbrac{k^{2}\slash\paren{\bZ\slash n\bZ}}^{\el}$. As such, one obtains an algebro-geometric description of the representation theory which can be fruitfully exploited in both directions. 
}
\bigbreak
It is hoped that one can extend the Azumaya splitting obtained in the previous section to formal neighbourhoods of fibers of the resolution map $\Psi:\ml\longrightarrow\Phi_{K^{\el}}^{-1}\paren{\eta^{l}}\slash K^{\el}$, and hence be able to combine the resulting abelian equivalence with the derived localization to obtain a description of module categories related to $\cU_{q,\eta}-\operatorname{mod}$ in terms of the geometry of $\ml$.\footnote{Such an extension to formal neighbourhoods of fibers or, alternatively, to derived fibers is necessary since the resolution map is not flat.} While such an equivalence would generalize that which is obtained using the methods described above for Cherednik algebras and resolutions of Kleinian singularities, we of course do not have as concrete a description of the global sections algebra or of the geometry in this more general setting. Obtaining a more explicit description of these objects, in order to gain a finer understanding of the derived equivalence, will be one focus of future work. \\

\addcontentsline{toc}{section}{References}
\bibliographystyle{plain}

\vspace{0.5pc}
{\sc Department of Mathematics and Statistics, Lancaster University, Lancaster, UK}

{\it Email address}: {\tt n.cooney@lancaster.ac.uk}
\end{document}